\documentclass[journal]{IEEEtran}
\usepackage{amsmath,amssymb}
\usepackage{graphicx}
\usepackage{cite}
\usepackage{hyperref}
\usepackage{epstopdf}
\usepackage{subfigure}

\newcommand{\comment}[1]{}

\comment{Kasemsak Uthaichana is with Department of Electrical
Engineering, Chiang Mai University, Chiang Mai, Thailand, (e-mail:
kasemsak@chiangmai.ac.th).  He is supported by grant Thailand Research Fund, Commission on Higher Education.\\
Sorin Bengea is with Control Systems Group, Research Center,
United Technologies, East Hartford, Connecticut, USA (email:sbengea@ieee.org)\\
Raymond DeCarlo and Steve Pekarek are with School of Electrical and
Computer Engineering, Purdue University, West Lafayette, Indiana, USA
(email: decarlo@ecn.purdue.edu, pekarek@purdue.edu)\\
Milo\v{s} \v{Z}efran is with Department of Electrical and Computer
Engineering, University of Illinois at Chicago, Chicago, Illinois, USA
(email: mzefran@uic.edu)}

\makeatletter
\newcommand{\submitted}[1]{\setbox\@tempboxa\vbox{\small \raggedright
    #1 \\ \hbox{}}
    \vspace*{-1.5 cm} \usebox\@tempboxa \\
    \vspace{-\ht\@tempboxa} \vspace{-\baselineskip} \vspace*{1.5 cm}}
\makeatother

\begin{document}

\title{\submitted{Preprint, Appeared in Journal of Nonlinear Systems
    and Applications, Vol. 2, No. 2, pp. 96--110, 2011. Copyright:
    Watam Press}
\Large\bf \uppercase{Hybrid Optimal Theory and Predictive Control for Power
Management in Hybrid Electric Vehicle}}
\author{Kasemsak Uthaichana\thanks{Kasemsak Uthaichana is with Department of Electrical
Engineering, Chiang Mai University, Chiang Mai, Thailand, (e-mail:
kasemsak@chiangmai.ac.th).  He is supported by grant Thailand
Research Fund, Commission on Higher Education.}, Raymond
DeCarlo\thanks{Raymond DeCarlo and Steve Pekarek are with School of
Electrical and Computer Engineering, Purdue University, West
Lafayette, Indiana, USA (email: decarlo@ecn.purdue.edu,
pekarek@purdue.edu)}, Sorin Bengea\thanks{Sorin Bengea is with
Control Systems Group, Research Center, United Technologies, East
Hartford, Connecticut, USA (email:sbengea@ieee.org)}, Milo\v{s}
\v{Z}efran\thanks{Milo\v{s} \v{Z}efran is with Department of
Electrical and Computer Engineering, University of Illinois at
Chicago, Chicago, Illinois, USA (email: mzefran@uic.edu)}, and Steve
Pekarek}
\date{}

\maketitle

{\footnotesize \noindent {\bf Abstract.}
  This paper presents a nonlinear-model based hybrid optimal control
  technique to compute a suboptimal power-split strategy for
  power/energy management in a parallel hybrid electric vehicle
  (PHEV).  The power-split strategy is obtained as model predictive
  control solution to the power management control problem (PMCP) of
  the PHEV, i.e., to decide upon the power distribution among the
  internal combustion engine, an electric drive, and other
  subsystems. A hierarchical control structure of the hybrid vehicle,
  i.e., supervisory level and local or subsystem level is assumed in
  this study.  The PMCP consists of a dynamical nonlinear model, and a
  performance index, both of which are formulated for power flows at
  the supervisory level.  The model is described as a bi-modal
  switched system, consistent with the operating mode of the electric
  ED.  The performance index prescribing the desired behavior
  penalizes vehicle tracking errors, fuel consumption, and frictional
  losses, as well as sustaining the battery state of charge (SOC).
  The power-split strategy is obtained by first creating the embedded
  optimal control problem (EOCP) from the original bi-modal switched
  system model with the performance index.  Direct collocation is
  applied to transform the problem into a nonlinear programming
  problem.  A nonlinear predictive control technique (NMPC) in
  conjunction with a sequential quadratic programming solver is used
  to compute suboptimal numerical solutions to the PMCP.  Methods for
  approximating the numerical solution to the EOCP with trajectories
  of the original bi-modal PHEV are also presented in this paper.  The
  usefulness of the approach is illustrated via simulation results on
  several case studies.\\
{\bf Keywords.} Hybrid Optimal Control, Nonlinear Model
Predictive Control, Hybrid Electric Vehicles, Power management,
Nonlinear Modeling}

\section{INTRODUCTION}

In a hybrid propulsion system, power distribution from two or more
energy sources/storages coordinate to deliver the performances demanded
by the drivers while considering fuel efficiency and operational
constraints.  In a parallel hybrid electric vehicle (PHEV), the power
demand can be delivered by the main power converter and/or the
energy-storage device.  Such energy storage devices could be batteries
with or without supercapacitors \cite{Erdinc2009, Pesaran2009}.
Examples of main power converters are internal combustion engines
(ICEs), fuel cells \cite{Barbir2005, Cao2009, Feroldi2009, Shen2008,
Springer1991}, etc.  In any case, as illustrated in \cite{Ao2007,
Lin2003, Pisu2007, Rizzoni1999}, the power distribution among the main
PHEV subsystems is computed at the supervisory level.  The model of the
PHEV at the supervisory level in this investigation is represented as a
bi-modal switched system, as opposed to models with higher number of
modes.

The description of the PMCP for constructing the model-based control
strategies consists of the PHEV dynamical model, and a performance
index (PI), both of which are formulated at the supervisory level.
Approaches to solve the PMCP in the literatures can be categorized
according to computational requirements as the real-time
implementable type, and the global optimal type.  The dynamic
programming (DP) approaches compute optimal solutions over the
driving cycles \cite{Ao2007, Lin2003, Scordia2005}.  The curse of
dimensionality of DP is well known.   Thanks to recent advances in
optimization, approximation approaches have been developed and
alleviate this problem \cite{Lee2008, Powell2009, Kim2010}.  Since
full knowledge of the driving cycles is still required, control
using DP is not real-time implementable.  Nevertheless, the results
can be used as benchmarks for comparing the degree of optimality
under replicated driving conditions.

Real-time implementable control strategies for the HEV, not optimal
over driving cycles, usually undergo fine-tuning on the actual vehicles
for desired performances under various assumptions and driving
conditions.  The list includes but is not limited to classical
instantaneous/static optimization, adaptive equivalent fuel consumption
minimization strategy (A-ECMS) \cite{Pisu2007}, simplified rule based ,
fuzzy logic based \cite{Erdinc2009, Schouten2003, Xiong2009}, and neural
network based \cite{Cao2009}.

The Nonlinear Model Predictive Control (NMPC) technique can provide
suboptimal solutions with respect to the PI over the predictive-window.
The degree of the optimality of the NMPC strategy is bounded by the
instantaneous and the global optimizations.  Note that the NMPC still
requires a few predictive partitions within the optimization window
(preview) of the driving profile, but
not as extensively as the dynamic programming approach.  The problem
underlying the NMPC strategy for PMCP in \cite{Kermani2008} is a mixed
integer optimization problem, e.g., \cite{Bemporad1999}, which is
computationally expensive.

In this study, the embedding technique in \cite{Bengea2005} is adopted to
formulate the PMCP as a (convex) embedded optimal control problem,
EOCP, from the original (non-convex) switched optimal control problem,
SOCP.  Hence, the degree of complexity for the embedded version of the
NMPC problem is lower.

Numerical methods for solving optimization problems include single
shooting, multiple-shooting and direct collocation \cite{Neuman1973,
Schaefer2007, VonStryk1993, Zefran1996}, etc.  Difficulties with the
single shooting method for the bi-modal PHEV is illustrated in
\cite{Uthaichana2005}.  Therein the necessary conditions are used to
solve for the optimal controls.  A superior version, called multiple
shooting method, is adopted to solve an optimization problem in
\cite{Diehl2006}.  Instead of dealing with adjoint equations as in
the multiple-shooting method, the direct collocation is adopted in
this investigation.  The embedding technique in conjunction with the
direct collocation method is used to transform the problem into an
NLP.  The numerical solution to the NLP is computed using sequential
quadratic programming (SQP) over a predictive window. More details
on other methodologies to obtain solutions to NMPC problems can be
found in \cite{DeHaan2006, Martinsen2004}.

The following describes the paper organization.  Section 2 summarizes
bi-modal switched model for the PHEV.  The performance index is
detailed in Section 3.  The PMCP is formulated as a multi-objective
embedded optimal control problem at the supervisory level in Section 4.
 Section 5 reviews the embedding technique and presents sufficient
conditions for existence of optimal solutions.  Section 6 describes the
numerical methodology.  Section 7 presents the hybrid optimal and NMPC
solutions for a sawtooth driving profile, and NMPC solutions for the
EPA highway and US06 supplemental FTP (EPA high-speed) driving
profiles.

\section{PHEV as a Bi-Modal Switched System}

The computation at the supervisory level is done based on the
presumption that the desired power level can be implemented at the
subsystem level.  The local closed loop controllers must track the
corresponding reference power demand, thereby decoupling the supervisory
and local level control problems.  Hence, the modeling at the
supervisory level should reflect the closed loop behaviors of the
subsystems.

\subsection{Summary of Hardware Descriptions}

\begin{figure}[h]
\centerline{\includegraphics[width=226pt]{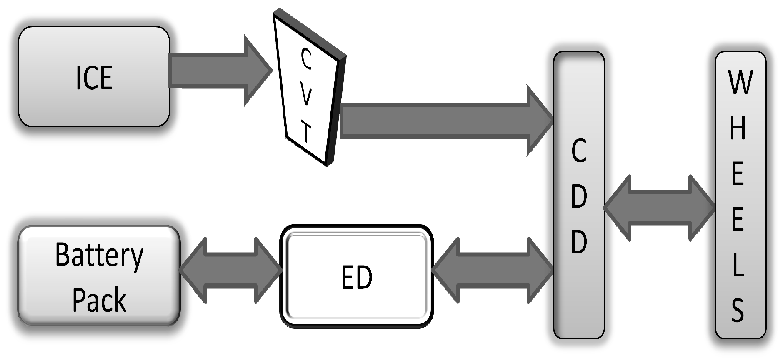}} \caption{Power flow diagram of PHEV in this study}
\label{fig1}
\end{figure}

The main power source is the 1.9 L ICE, coupled to the driveshaft
through a continuously variable transmission (CVT) and clutch in
the post-transmission configuration.  Thirty 13Ah 12V, lead-acid
batteries in series are interfaced with a 30 kW induction ED.  Hotel
loads are handled by a traditional engine-based charging system. The
coupling device and differential (CDD) acts as a summing junction for
redirecting the power flow among the ICE-CVT, battery-ED, and the
wheels. Figure \ref{fig1} illustrates power distributions among the
main subsystems.

\subsection{Modes of Operation}

In \cite{Lin2003}, five modes of operation describe the essential
behavior of the PHEV.  Since the PMCP complexity increases
exponentially with increased numbers of modes \cite{Floudas1995,
Nemhauser1988}, an effort is put forth in this investigation to reduce the
number of modes at the supervisory level.

\begin{figure}[h]
  \centerline{\includegraphics[width=226pt]{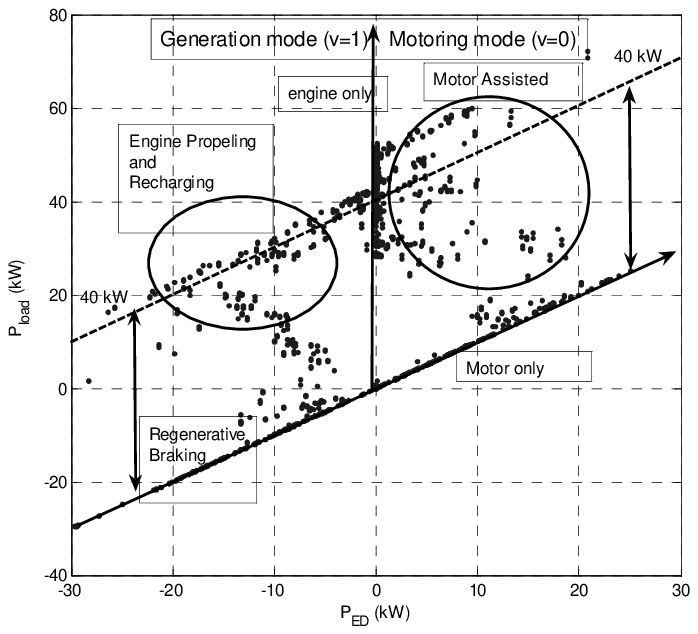}}
\caption{Power Flow strategy obtained in this study using NMPC tracking EPA highway profile showing the concept of mode
reduction from five to two}
\label{fig2}
\end{figure}

Through careful consideration of the dynamics in each mode of
operation, the essential behavior can be approximated using only two
modes.  As a preview, the mode reduction concept from five to two is
illustrated via the numerical results obtained in this study in
Fig. \ref{fig2}.  In Fig. \ref{fig2},
${{P}_{ED}}+{{P}_{ICE}}={{P}_{load}}$, i.e., the sum of ED power, and
ICE power is delivered to the load (the planetary efficiency is
ignored for now).

The details of the mode reduction concept, when the HEV is operating,
can be described as follows:

\begin{itemize}
\item For $v=0$  (${{P}_{ED}}\ge 0$): the engine-only mode
  $({{P}_{ICE}}>0,{{P}_{ED}}=0)$ corresponds to zero power flow from
  the ED which can be achieved by a zero-value of the ED control
  variable; the motor only mode (${{P}_{ICE}}=0,{{P}_{ED}}>0$)
  corresponds to no power flow from the ICE that is also achievable by
  a zero-value of the engine control variable; in motor assisted mode
  (${{P}_{ICE}}>0,{{P}_{ED}}>0$) both ICE power and ED power are
  strictly positive to the wheels achievable by non-zero control
  variables of the ED and ICE.
\item For $v=1$  (${{P}_{ED}}<0$): the regenerative-braking mode
  (${{P}_{load}}<0,{{P}_{ED}}<0$) corresponds to a reverse of the ED
  power flow to charge the battery with ICE power at zero (or nearly
  so); and the engine-charging-battery mode
  (${{P}_{ICE}}>0,{{P}_{ED}}<0$) corresponds to the case when the ED
  operates as a generator with positive ICE power flow to the ED and
  possibly to the wheels.
\end{itemize}

It can be seen that the two modes of operation at the supervisory
level coincide with the modes of the ED denoted as  $v=0$ (motoring)
and  $v=1$  (generating).  The summary of the bi-modal switched
system describing the essential dynamics of the power flow at the
supervisory level is given next.

\subsection{HEV State Model Overview}

The summary of the modeling equations, detailed in
\cite{Uthaichana2006, Uthaichana2008}, are summarized in this
subsection.  The input-output power relationships of subsystems
consist of algebraic and dynamical equations.  The input-output
relationship is considered algebraic when the its internal power
flow dynamics are much faster than others.  The essential dynamical
state is
\begin{equation}
x(t)\triangleq {{\left[ {{P}_{ICE}},SOC,V \right]}^{T}},
\end{equation}
where ${{P}_{ICE}}$ is the ICE power, $SOC$ is the battery state-of-charge, and $V$ is the longitudinal vehicle's
velocity.  The mode-dependent nonlinear state equation for the PHEV in this study is:
\begin{equation}
\dot{x}(t)={{f}_{v(t)}}(x(t),{{u}_{v(t)}}(t))
\label{ZEqnNum543948}
\end{equation}
where ${{f}_{v}}(\cdot )$ denotes the dynamics when motoring, $v=0$, or when generating, $v=1$.  The modulating controls
in modes-0/1 are
\begin{equation}
\begin{array}{c}
u_{0/1}(t)={{\left[ u_{ICE}(t),u_{FR}(t),u_{EM/GEN}(t) \right]}^{\text{T}}} \\
\ \ \ \ \ \in \Omega \subset {{\mathbb{R}}^{3}} \\
\end{array}
\end{equation}
where (i) the compact and convex set,
\begin{equation}
\Omega =[0,1] \times [0,1] \times [0,1]
\label{ZEqnNum444157}
\end{equation}
(ii) ${{u}_{ICE}}(t)\in [0,1]$ modulates the maximum available ICE
power; (iii) ${{u}_{FR}}(t)$ modulates the maximum frictional braking;
and (iv) ${{u}_{EM/GEN}}(t)$ $\in [0,1]$ modulates the maximum
available ED power in the mode-0 and mode-1, respectively.

The motivation for this control structure is four-fold.  First, the
model is scalable in terms of numbers of power sources (ICE's or ED's)
so that the corresponding increase in the number of operating modes
leads to only a polynomial increase in complexity for numerical
optimization methods \cite{Gill1981}.  Second, the model has a form
compatible with hybrid optimal control theory.  Third and more
critically, the fact that the controls take values in a convex compact
set  $\Omega $  makes the PMCP amenable to hybrid optimization techniques.
Fourth, any optimization algorithm searches for the optimal controls
and switching function $v(t)$ in a hypercube as opposed to a
(non-convex) state and time-dependent region.

\subsubsection{State Equation for the ICE}

The variable  ${{P}_{ICE}}$  denotes the unidirectional instantaneous
ICE power flow, quantified at the flywheel and including losses due to
parasitic loads.  The ICE dynamical equation is given by equation
\eqref{ZEqnNum262617}, i.e.,
\begin{equation}
\begin{array}{c}
{{{\dot{P}}}_{ICE}}=-\frac{1}{{{\tau }_{ICE}}}{{P}_{ICE}}+\frac{1}{{{\tau }_{ICE}}}P_{ICE}^{\max }\left(
{{\omega }_{ICE}} \right) \\
\ \ \ \ \ \ \ \ \ \ \ \ \ \ \ \ \ \cdot eng\left( {{\omega }_{ICE}} \right)\cdot {{u}_{ICE}}(t) \\
\label{ZEqnNum262617}
\end{array}
\end{equation}
where ${{\tau }_{ICE}}$ is the nominal engine power delivery delay averaging the effect of the firing delay, smoke limit
map, crankshaft speed, fueling shot mode, etc.  Further, this also ensures that the command handed down by the supervisory
controller can be followed. $P_{ICE}^{\max }$ is an  ${{\omega }_{ICE}}$ -dependent maximum available ICE power;
${{\omega }_{ICE}}$  is the CVT controller-selected engine speed using the strategy modified slightly from the
speed-envelope for a non-hybrid ICE in \cite{Pfiffner2001}. Specifically,
\begin{equation}
{{\omega }_{ICE}}=\left( 1-p \right)\omega _{ICE}^{\min }(V)+p\omega
_{ICE}^{\max }(V)
\label{ZEqnNum809968}
\end{equation}
where  $p\in [0,1]$  modulates the speed curve according to the ICE
power level;  $\omega _{ICE}^{\min }(V)$  and  $\omega _{ICE}^{\max
}(V)$  are the minimum and the maximum allowable speeds at each
vehicle's velocity and illustrated in Fig. 3.  For better driveability,
ICE-CVT's transition from non-engaged to engaged is not allowed when
the vehicle's velocity is too low for jerk reduction.  Note that in
this study, the capitalized superscript indicates constant whereas
lower-case superscripted means parameter dependent.

\begin{figure}[h]
\centerline{\includegraphics[width=230pt]{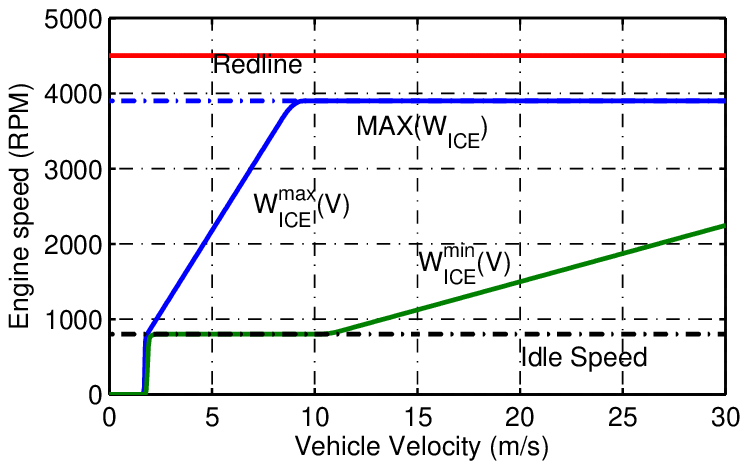}} \caption{Maps of minimum and maximum engine speeds for each vehicle's
velocity}
\label{fig 3}
\end{figure}

According to ICE dynamical equation \eqref{ZEqnNum262617}, and the constraints on the input, \eqref{ZEqnNum444157}, it can
be shown that  ${{P}_{ICE}}$  belongs to the compact and convex set:
\begin{equation}
{{P}_{ICE}}\in [P_{ICE}^{MIN},P_{ICE}^{MAX}]\subset
\mathbb{R}
\label{ZEqnNum990052}
\end{equation}
The range of ICE power \eqref{ZEqnNum990052} is not an additional state constraint, but rather is a direct result from the
aforementioned relationships.  This information is useful during the computation of the solution as the numerical
search-space is smaller.

\subsubsection{State Equation for Battery Operation}

For a relatively new battery, the normalized energy and the SOC are
equivalent \cite{Coleman2008, Rao2003}.  The state equation describing
the SOC dynamics is derived based on the conservation of power and
energy.  The parameters are computed to fit the battery data.  The
\textit{partial} linearization about the mode-dependent nominal battery
operating power, $P_{bat,nom}^{v}$, describing the SOC is given by:
\begin{equation}
\begin{array}{l}
\overset{\centerdot }{\mathop{SOC}}\,(t)=\frac{{{d}_{3,v}}}{W_{bat}^{MAX}}{{\left( P_{bat,nom}^{v} \right)}^{2}} \\
\ \ \ -\left[ \ln \left( {{d}_{2,v}} +{{d}_{1,v}}SOC(t) \right)+2{{d}_{3,v}}P_{bat,nom}^{v} \right. \\
\ \ \ \left. +{{d}_{4,v}} \right]\frac{{{P}_{bat}}}{W_{bat}^{MAX}}
\label{ZEqnNum706757}
\end{array}
\end{equation}
The validation result of this model against a variety of battery data appear in
\cite{Agarwal2009}.  In equation \eqref{ZEqnNum706757}, $W_{bat}^{MAX}$  is the rated maximum battery energy,
${{P}_{bat}}$  is the battery power either drawn by (positive for $v=0$) or provided by the ED (negative for $v=1$) and is
implicitly controlled by ${{u}_{EM/GEN}}$.  Specifically,
\begin{equation}{{P}_{bat}}=\left\{ \begin{matrix}
P_{ED,in}^{0}, & v=0  \\
-P_{ED,out}^{1}, & v=1  \\
\end{matrix} \right.
\end{equation}
${{d}_{k,v}}$, $k=1,...,4$ are the appropriate coefficients.  The consideration of the recovery, cycling, and aging
effects are beyond the scope of this investigation.  This formulation makes equation \eqref{ZEqnNum706757} scalable to a
variety of battery storage capacities and types.

\subsubsection{State Equation for Vehicle Motion}

The conventional longitudinal vehicle's velocity is described, not in
terms of torques, but rather in terms of the acting power flow as:
\begin{equation}
\begin{array}{l}
\dot{V}=-\left[ \frac{{{k}_{v1}}}{{{m}_{c}}}{{V}^{2}}+{{k}_{v2}}\cos \left( \alpha(t) \right) \right]sgn\left( V \right) \\
-g\sin \left( \alpha (t) \right)+\frac{1000}{{{m}_{c}}\left( V+{{\varepsilon }_{V}} \right)}\left[ P_{CDD,wh}^{v}-{{P}_{FR}} \right] \\
\label{ZEqnNum493755}
\end{array}
\end{equation}
In equation \eqref{ZEqnNum493755},  ${{\varepsilon }_{V}}$ is a regularization term; ${{m}_{c}}$  is vehicle mass;
$\frac{{{k}_{v1}}}{{{m}_{c}}}{{V}^{2}}$ is normalized aerodynamic drag; ${{k}_{v2}}\cos (\alpha (t))$  is the rolling
resistance; $\alpha (t)$is the time-varying angle of road inclination; $P_{CDD,wh}^{v}$  is the power delivered from ($\ge
0$) and to ($<0$) the CDD. Finally,
\begin{equation}
{{P}_{FR}}=P_{FR}^{\max }(V){{u}_{FR}}(t)
\end{equation}
is the frictional braking power. As a result from equation \eqref{ZEqnNum493755}, the vehicle's velocity is also in an
invariant set,
\begin{equation}
V\in \left[ {{V}^{MIN}},{{V}^{MAX}} \right]\subset R.
\end{equation}

\subsubsection{Mode Dependent ED Modeling Equations}

The derivation of the ED algebraic input-output power flow equations
for both modes can be found in \cite{Pekarek2005}.  The ED in this
study, operated under a maximum torque/amp (MTA) control strategy, can
be represented at the supervisory level as
\begin{equation}
P_{ED}^{v}=\eta _{ED}^{v}({{\omega
}_{ED}})P_{ED,in}^{v}
\label{ZEqnNum174072}
\end{equation}
Each term in equation \eqref{ZEqnNum174072} is mode dependent. The ED output power is denoted  $P_{ED}^{v}$ , the
efficiency $\eta _{ED}^{v}({{\omega }_{ED}})$  strongly depends on the choice of closed-loop control, a phenomena largely
underweighted in the HEV literature; ${{\omega }_{ED}}=\beta \cdot V$  is the ED rotor speed;  $\beta $  is a positive
constant. The ED input power in modes 0 and 1 is
\begin{equation}
P_{ED,in}^{0}=P_{ED,in}^{max}\left( {{\omega }_{ED}} \right)\cdot {{u}_{EM}}(t)
\end{equation}
\begin{equation}
P_{ED,in}^{1}=P_{ED,in}^{max}\left( {{\omega }_{ED}} \right)\cdot {{u}_{GEN}}(t)
\end{equation}
where $P_{ED,in}^{max}\left( {{\omega }_{ED}} \right)$  is the speed dependent ED maximum input power modulated by the
control ${{u}_{EM}}(t)$  in mode-0, and ${{u}_{GEN}}(t)$ in mode-1.

\subsubsection{CVT and mode-dependent CDD Power Flow Equations}

No power response lag between the input and output CVT powers is
assumed at the supervisory level, leading to the algebraic equation
\begin{equation}
{{P}_{cvt,out}}(t)={{\eta }_{cvt}}{{P}_{cvt,in}}
\end{equation}
where ${{\eta }_{cvt}}$ is the CVT efficiency; ${{P}_{cvt,in}}(t)={{P}_{ICE}}$ is the CVT input power; and the output
power is delivered to the CDD, i.e., ${{P}_{cvt,out}}(t)={{P}_{CDD,cvt}}$.

The CDD's input/output power flows are given by
\begin{equation}
P_{CDD,wh}^{0}(t)={{\eta }_{cdd1}}{{P}_{CDD,cvt}}+{{\eta
}_{cdd2}}P_{CDD,ED}^{0}
\end{equation}
and
\begin{equation}
P_{CDD,ED}^{1}(t)={{\eta }_{cdd2}}{{P}_{CDD,cvt}}-{{\eta
}_{cdd2}}P_{CDD,wh}^{1}
\end{equation}
(i) ${{\eta }_{cdd1}}$, and ${{\eta }_{cdd2}}$ are the appropriate
power transfer efficiency among the ED, CVT and wheels; (ii)
$P_{CDD,ED}^{0}=P_{ED}^{0}$  is the propulsion power coming directly
from the output of the ED in mode-0; in mode-1,
$P_{CDD,ED}^{1}=P_{ED,in}^{1}$  is an output power port of the CDD
providing mechanical power to the input of the ED (generator); (iii)
in mode-0,  $P_{CDD,wh}^{0}\ge 0$.  However, in mode-1,
$P_{CDD,wh}^{1}(t)$  can be either positive or negative. Note that
$P_{CDD,wh}^{0}$ is represented as ${{P}_{load}}$ in the mode
reduction concept in Section 2.2.

\section{PERFORMANCE INDEX}

To incorporate the desired behaviors of the HEV operation, we consider the optimization functional for each mode, as
follows:

\begin{equation}
\begin{array}{rcl}
{{J}_{v}}({{x}_{0}},u,[{{t}_{0}},{{t}_{f}}]) & = &
g({{t}_{0}},{{x}_{0}},{{t}_{f}},{{x}_{f}})\\
& & +\int_{{{t}_{0}}}^{{{t}_{f}}}{{{L}_{v}}(t,x,u)dt}
\end{array}
\label{ZEqnNum287788}
\end{equation}

The mode-dependent integrand  ${{L}_{v}}(t,x,u)$  depends on the
optimization objectives, such as minimizing only fuel consumption as
in \cite{Anatone2005}, or a combination of fuel consumption and
emissions as in \cite{Jeon2002, Koot2005, Lin2003}.  In this
research the PI consists of terms that are consistent with the power
flow management framework and have meaningful physical
interpretations.  The integral quadratic PI that uses the same
integrand for both modes of operation, i.e.,
${{L}_{0}}(t,x,u)={{L}_{1}}(t,x,u)$ is adopted in this study.  The
integrand for both modes is
\begin{equation}
\begin{array}{rcl}
{{L}_{v}} & = & {{C}_{V}}{{\left( V-{{V}^{des}}(t) \right)}^{2}} \\
& & +{{C}_{ICE}}{{\left(
\frac{{{P}_{ICE}}}{{{\eta}_{ICE}}(\centerdot )}
\right)}^{2}}+{{C}_{FR}}{{\left( {{P}_{FR}} \right)}^{2}}
\end{array}
\end{equation}

The integrand penalizes the velocity tracking error,
${{C}_{V}}{{\left( V-{{V}^{des}}(t) \right)}^{2}}$ , the frictional
braking power,  ${{C}_{FR}}{{\left( {{P}_{FR}} \right)}^{2}}$ , and
the fuel usage.  The fuel usage is approximated by ICE power usage
divided by fuel conversion efficiency \cite{Heywood1988}, i.e.,
${{C}_{ICE}}{{\left( \frac{{{P}_{ICE}}}{{{\eta }_{ICE}}(\centerdot
)} \right)}^{2}}={{C}_{ICE}}{{({{P}_{fuel}})}^{2}}$  where  ${{\eta
}_{ICE}}\left( {{P}_{ICE}},V \right)$  is the ICE efficiency that
depends on the ICE power-and-speed.  Fig. 4 depicts the efficiency
map of the ICE superimposed with the iso-efficiency curves.
\begin{figure}[h]
\centerline{\includegraphics[width=229pt]{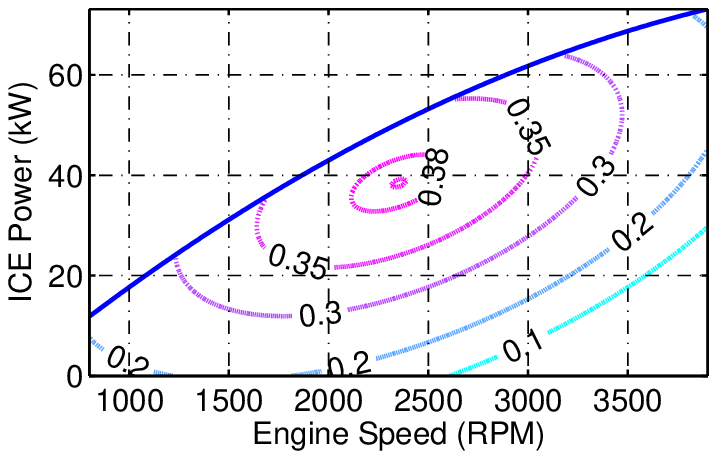}} \caption{ICE
power-vs.-engine speed superimposed with iso-efficiency curves}
\label{Fig 4}
\end{figure}

The penalty on the variation in the boundary conditions,
$g({{t}_{0}},{{x}_{0}},{{t}_{f}},{{x}_{f}})$, in PI
(\ref{ZEqnNum287788}) is taken as ${{C}_{bat}}(\centerdot ){{\left(
SOC({{t}_{f}})-SO{{C}^{NOM}} \right)}^{2}}$.  This choice of the
penalty pushes the SOC at ${{t}_{f}}$  toward the nominal level,
$SO{{C}^{NOM}}$.  It is also desirable to operate the SOC in a
predefined range to prolong battery lifetime.  It will be shown in
the simulation results that this choice of penalty on the battery
SOC can be used to encourage
\begin{equation}
SOC\in \left[ SO{{C}^{MIN}},SO{{C}^{MAX}} \right]\subset
\mathbb{R}
\label{ZEqnNum384721}
\end{equation}

Note that if the constraints (\ref{ZEqnNum384721}) are violated, the
penalty term on the SOC must be more stringent.  This SOC strategy
is intended to enforce a charge-sustaining operation.

A more elaborate PI accounting for drivetrain losses in each mode
has the form
\begin{equation}
\begin{array}{l}
{{L}_{v}}= C_{V}^{v}\left( V-{{V}^{des}}
\right)+C_{ICE}^{v}{{\left( {{P}_{fuel}} \right)}^{2}} \\
\ \ +C_{cvt}^{v}{{\left( {{P}_{cvt,loss}} \right)}^{2}}
+C_{CDD}^{v}{{\left( P_{C,loss}^{v}\right)}^{2}} \\
\ \ +C_{ED}^{v}{{\left( P_{ED,loss}^{v} \right)}^{2}}
+C_{bat}^{v}{{\left( P_{bat,loss}^{v} \right)}^{2}} \\
\ \ +C_{FR}^{v}{{\left( {{P}_{FR}} \right)}^{2}}
\end{array}
\end{equation}
where the additional power loss terms are CVT losses, CDD loses, ED
losses, and battery losses, whose identity should be clear from the
notation.  More details on the generalized PI can be found in
\cite{Uthaichana2006}.

\section{POWER MANAGEMENT CONTROL PROBLEM AND THE EOCP}

For the switched optimal control problem (SOCP), the modal switching
function $v(t)$  belongs to a discrete set  $\{0,1\}$,  $v(t)\in
\{0,1\}$.  In contrast, for the EOCP the modal switching function
$v(t)$  takes values in the closed interval $[0,1]$, a continuum of
possible values.  The enlargement of  $v(t)\in \{0,1\}$  to $v(t)\in
[0,1]$ constitutes an embedding of the SOCP into a larger family of
continuously parameterized problems.  This embedding converts a
non-convex SOCP into a convex EOCP.  As per \cite{Bengea2005} the
SOCP can almost always be solved by first solving the EOCP and any
solution of the EOCP can be approximated to any degree of precision
by some solution of the switched state model \eqref{ZEqnNum543948}.
Further, in this study, projection techniques are also presented as
alternatives for approximating the EOCP solution by an SOCP
trajectory.

\subsection{Specification of the embedded optimal control problem}

The embedding requires creating a convex combination of the vector
fields of the switched state model according to the equation,
\begin{equation}
\begin{array}{rcl}
\dot{x}(t) & \triangleq &{{f}_{E}}\left(
x(t),{{u}_{0}}(t),{{u}_{1}}(t),v(t) \right) \\
& = &\left[ 1-v(t) \right]{{f}_{0}}\left( x(t),{{u}_{0}}(t)
\right) \\
& & +v(t){{f}_{1}}\left( x(t),{{u}_{1}}(t) \right) \\
\label{ZEqnNum853462}
\end{array}
\end{equation}
where ${{u}_{i}}(t)\in \Omega $, $i=0,1$.  Clearly if $v(t)=0$,
${{f}_{E}}$ reduces to the 0-mode vector field and similarly
for $v(t)=1$.

The performance index (PI) of the EOCP results from a similar convex
embedding of the PIs associated with each mode of operation of the
SOCP:

\begin{equation}
\begin{array}{l}
{{J}_{E}}({{x}_{0}},{{u}_{0}},{{u}_{1}},v,[{{t}_{0}},{{t}_{f}}])=g({{t}_{0}},{{x}_{0}},{{t}_{f}},{{x}_{f}}) \\
\ \ \
+\int_{{{t}_{0}}}^{{{t}_{f}}}{{{L}_{E}}(t,x,{{u}_{0}},{{u}_{1}},v)dt}
\\
\ \ \ =g({{t}_{0}},{{x}_{0}},{{t}_{f}},{{x}_{f}}) \\
\ \ \ +\int_{{{t}_{0}}}^{t_f}{\left[ \left( 1-v(t)
\right){{L}_{0}}(t,x,{{u}_{0}})\right.} \\
\ \ \ \left. +v(t){{L}_{1}}(t,x,{{u}_{1}}) \right]dt
\end{array}
\label{ZEqnNum693203}
\end{equation}
with ${{L}_{i}}(t,x,{{u}_{i}})$, $i=0,1$, denoting the
convex-in-${{u}_{i}}$ integrands of the PI.  When  $v(t)\in
\{0,1\}$, the minimization of (\ref{ZEqnNum693203}) subject to
(\ref{ZEqnNum853462}) defines the SOCP, while when  $v(t)\in [0,1]$,
the minimization of (\ref{ZEqnNum693203}) subject to
(\ref{ZEqnNum853462}) constitutes the EOCP.  Formally the EOCP (the
structure for solving the PMCP) becomes:
\begin{equation}
\underset{{{u}_{0}},{{u}_{1}},v}{\mathop{\text{min}}}\,\text{
}{{J}_{E}}({{x}_{0}},{{u}_{0}},{{u}_{1}},v,[{{t}_{0}},{{t}_{f}}])
\end{equation}
with ${{J}_{E}}(\cdot )$ given by (\ref{ZEqnNum693203}), subject to
\begin{equation}
\dot{x}(t)={{f}_{E}}(x(t),{{u}_{0}}(t),{{u}_{1}}(t),v(t))
\end{equation}
with  ${{f}_{E}}$  given in (23), $v(t)\in [0,1]$, and
${{u}_{0}}$, ${{u}_{1}}\in \Omega $.

\subsection{Relationships between EOCP and SOCP}

If the EOCP has a bang-bang type solution (wherein  $v(t)$  only
takes values in $\{0,1\}$ ) then clearly it is also a solution to
the original SOCP.  Further it can be shown (Corollary 2 in
\cite{Bengea2005}) that the set of trajectories of the switched
system (equation (\ref{ZEqnNum853462}) with  $v(t)\in \{0,1\}$) is
dense in the set of trajectories of the embedded system (equation
(23) with $v(t)\in [0,1]$ ).  Thus when/if the EOCP does not have a
bang-bang type solution (wherein  $v(t)\in (0,1)$  for some non-zero
measure sets of time) then the EOCP solution can be approximated by
a trajectory of the switched system to any desired degree of
precision.  These relationships between the SOCP and EOCP motivate
and justify the effort in determining SOCP solutions by solving the
EOCP. Additional relationships between SOCP and EOCP can be found in
\cite{Bengea2005}.

\subsection{Approximation to Singularities in EOCP}

This subsection describes approximation techniques when the control
$v(t)$ obtained via the EOCP is not bang-bang.  When, $v(t)\in (0,1)$,
i.e., $v(t)$ takes on fractional values, over an interval ${{t}_{1}}\le
t\le {{t}_{2}}$, it would suggest that for the HEV the ED operate
simultaneously in both modes for this time interval, an impossibility.
In another words, when $v(t)\in (0,1)$, the SOCP does not have a
solution, but epsilon-approximating solutions to the EOCP can be constructed as
follows.  Given a desired error of approximation, $\varepsilon $, one
can construct subintervals
${{t}_{1}}<{{T}_{1}}<{{T}_{2}}....<{{t}_{2}}$
such that $| {{T}_{i+1}}-{{T}_{i}}|<\delta $,
where delta is generated based on $\varepsilon $, vector fields,
${{f}_{i}}$, and cost integrand, ${{L}_{i}}$.  In the case when the
switching interval--constrained by the embedded controller loop time,
actuator bandwidth, etc--is larger than $\delta $, one would have to
increase the approximating error, $\varepsilon $, and re-construct the
intervals.  The approximating error will need to be sufficiently large
to accommodate the constraint $\delta >{{T}_{\min
}}$, where ${{T}_{\min }}$ is the minimum switching period.

The construction of the switching subintervals in the case of complex
vector fields, such as the case for the HEV model, can be alleviated by
considering empirical based switching intervals such as described
below.

One approach to empirical switching is to average the fractional
values of $v(t)$ over ${{t}_{1}}\le t\le {{t}_{2}}$ and the average
value, denoted $\bar{v}$, over ${{t}_{1}}\le t\le {{t}_{2}}$, can be
interpreted as a duty cycle, or a pulse width modulation (PWM)
control.  So there exists a time $t'$ such that for ${{t}_{1}}\le
t<t'$, the system is in mode 0 and for $t'\le
t<{{t}_{2}}$ the system is in mode 1 so that the
average over the whole interval is
$\bar{v}=\frac{{{t}_{2}}-t'}{{{t}_{2}}-{{t}_{1}}}$.  Thus, a PWM or
switched approximation to the embedded $v(t)$ is made.

The previously computed ${{u}_{i}}(t)$ are associated with the
embedded solution $v(t)$, not the new approximation.  One possibility
is to simply use these values for each associated subinterval.  A
second possibility is to set $v(t)$ equal to its PWM approximation and
then find the optimal ${{u}_{0}}$ and ${{u}_{1}}$ associated with this
choice.  Switching can be minimized by beginning the duty cycle for the
next interval in the ending mode of the prior interval.

A third possibility is as follows: let ${{T}_{\min }}$ be the smallest
switching interval of time.  For each time unit, ${{t}_{1}}\le t\le
{{t}_{2}}={{t}_{1}}+{{T}_{\min }}$, one can project the fractional
value of $v(t)$ onto the set $\{0,1\}$ according to the formula:

\begin{equation}
{{\| \overline{(1-v(t))\cdot {{u}_{0}}(t)} \|}_{2}}\left\{
\begin{matrix}
\ge {{\overline{\| v(t)\cdot {{u}_{1}}(t) \|}}_{2}} \Rightarrow  v(t)=0  \\
<{{\| \overline{v(t)\cdot {{u}_{1}}(t)} \|}_{2}} \Rightarrow v(t)=1
\end{matrix} \right.
\end{equation}
where over-bars denote averages over the interval ${{t}_{1}}\le t\le
{{t}_{2}}={{t}_{1}}+{{T}_{\min }}$.  As before one can either use
the previously calculated values of ${{u}_{i}}$ or resolve the
optimization with $v(t)$ fixed at the desired mode.  For the
simulation studies of this work, equation (27)  was used to fix the
bang-bang solution for $v(t)$ and then the optimization was resolved
for the best pair of ${{u}_{i}}(t)$ given the fixed mode sequence.

\subsection{Embedded PI for PMCP}

As mentioned earlier, the integrand and the penalty on the boundary
conditions in both modes are the same.  The embedded PI for the PMCP
is obtained by substituting appropriate terms in the PI
(\ref{ZEqnNum693203}), i.e.,
\begin{equation}
\begin{array}{rcl}
{{J}_{E}} & = & {{C}_{bat}}(\centerdot ){{\left( SOC({{t}_{f}})-SO{{C}^{NOM}} \right)}^{2}} \\
& & +\int_{{{t}_{o}}}^{{{t}_{f}}}{\left({{C}_{V}}{{\left( V-{{V}^{des}}(t) \right)}^{2}} \right.}\\
& & +{{C}_{ICE}}{{\left( \frac{{{P}_{ICE}}}{{{\eta }_{ICE}}(\centerdot )} \right)}^{2}} \\
& & \left. +{{C}_{FR}}{{\left( P_{FR}^{MAX}(\centerdot){{u}_{FR}}(t)
\right)}^{2}} \right)dt
\end{array}
\label{ZEqnNum948398}
\end{equation}
where the physical meaning of each term is given in Section 3.

\section{SUMMARY ON THEORETICAL FOUNDATIONS}

When the discrete input $v\in \{0,1\}$  presents, it renders, in
general, the SOCP non-convex.  For a variety of assumptions on
system vector fields ${{f}_{v}}$, an SOCP performance index, and
mode-switching penalties and constraints, several approaches have
been employed in the literature for characterizing and computing
SOCP solutions, consisting of: searches over or assumptions on mode
sequences and switching instants, after which one computes the
continuous control values and the cost to compare the different
scenarios.  These approaches do not allow the switching function to
be chosen in concert with the continuous time control as is the case
with the embedded approach.

Discussing neither sufficient conditions for optimality nor account
for the singular solution scenarios, Riedinger et al. (1999) applies
directly the Maximum Principle to the SOCP.  For a larger class of
systems, and with a cost that depends on the mode sequence, Sussmann
(1999) derives necessary conditions for optimality via a generalized
Maximum Principle.  Other approaches include pre-assigned switching
sequence method (for a limited class of problems) in \cite{Giua2001},
and a hybrid Bellman inequality approach in \cite{Hedlund1999}.  Mixed
integer programming (MIP) approaches have also been employed to find
optimal solutions \cite{Bemporad1999}.  Solving the problem using MIP
methods, however, is non-deterministic polynomial-time hard (NP-hard);
indeed the scalability of this technique is problematic \cite{Wei2007}.

The nonconvexity of the problem and the inapplicability of the
mentioned existing techniques\textemdash{}too general and impractical,
or very specific results, or insufficient characterization of
solutions\textemdash{} to the SOCP has led to the development of the
parameterized family of problems, the EOCP, set forth in the previous
section.

\subsection{EOCP: Sufficient Existence Conditions}

This section summarizes the main sufficient conditions for EOCP
solutions.  Sufficient conditions for optimality are [Theorem 9, in
\cite{Bengea2005}]:
\begin{itemize}
\item[(i)] the admissible pair set (control, trajectory) is nonempty;
\item[(ii)] the points $\left( t,x(t) \right)$ are included in a compact
set for all$t\in [{{t}_{0}},{{t}_{f}}]$;
\item[(iii)] the terminal set is compact;
\item[(iv)] the input constraint set is compact and convex;
\item[(v)] the vector fields ${{f}_{0}}$ and ${{f}_{1}}$ are linear in their
(control) inputs ${{u}_{0}}$, and ${{u}_{1}}$, respectively i.e.,
\begin{itemize}
\item[(S1)] ${{f}_{0}}(t,x,{{u}_{0}})={{A}_{0}}(t,x)+{{B}_{0}}(t,x){{u}_{0}}$
\item[(S2)] ${{f}_{1}}(t,x,{{u}_{1}})={{A}_{1}}(t,x)+{{B}_{1}}(t,x){{u}_{1}}$
\end{itemize}
\item[(vi)] for each$\left( t,x(t) \right)$, the integrands of the penalty
functions, ${{L}_{0}}(t,x,{{u}_{0}})$ and
${{L}_{1}}(t,x,{{u}_{1}})$, are convex functions of ${{u}_{0}}$, and
${{u}_{1}}$, respectively.
\end{itemize}

Based on the assumptions made on the input constraint set and on the
vector fields ${{f}_{0}}$ and ${{f}_{1}}$, one can conclude that
conditions (i), (ii), and (iv) are met.  Further, a sufficiently
large compact set can be substituted for the terminal set, meeting
condition (iii).  Condition (v) is also met as it can be observed
based on the modeling equations from Sections 2.  Specifically, the power terms
that depend on the continuous control inputs, are
factored into the product of a control input and a term that depends
on the state, $x(t)$. Utilizing the form of these power terms, and
the forms of ${{L}_{0}}(t,x,{{u}_{0}})$ and
${{L}_{1}}(t,x,{{u}_{1}})$ one concludes that condition (vi) is also
met.  Hence the EOCP has a solution.

The above sufficient conditions only guarantee the existence of the
EOCP's solutions, but do not provide a solution methodology.  In
conjunction with the SOCP-EOCP relationships mentioned above, the
necessary conditions obtained by direct application of the Maximum
Principle \cite{Pekarek2005} provide a method for obtaining at least
suboptimal solutions of the SOCP.   By using this approach, the
optimization problem is transformed into a two-point boundary value
problem on the state and adjoint equations.  The single shooting
method is applied to compute the numerical solution in
\cite{Uthaichana2005}, and the solution is very sensitive with
respect to the co-state initial condition.  The multiple shooting
method can be applied to reduce the sensitivity issue.  This paper
takes an alternate approach for computing numerical solutions to the
EOCP, i.e., via the direct collocation method, described in the next
section.

\section{NUMERICAL TECHNIQUE AND NMPC}

This section briefly describes the direct collocation method and the
nonlinear model predictive control (NMPC) strategy.  Both are used in
conjunction to formulate the PMCP as a nonlinear programming problem
(NLP).

\subsection{Discretization via Direct Collocation}

Given the PI (\ref{ZEqnNum948398}) and the state equation and
constraints of equations (\ref{ZEqnNum853462}), one discretizes
these equations using the direct collocation method.  The
discretization of the PI uses a variation of the trapezoidal rule
and constraint equations use the mid-point rule, respectively.
These discretized equations convert the EOCP into a finite
dimensional NLP where states and inputs are treated as unknown
variables.  The direct collocation technique consists of several
steps that have two main stages: (i) time discretization, and state
and input function approximations by a finite number of polynomial
basis functions; (ii) approximation of the continuous state dynamics
and cost index integrand by discrete-state and
discrete-input-dependent counterparts.

Without going through a lengthy derivation, the continuous time
interval $[{{t}_{0}},\,\,{{t}_{f}}]$ is discretized into a sequence
of points ${{t}_{0}}<{{t}_{1}}<{{t}_{2}}<\cdot \cdot \cdot
<{{t}_{N-1}}<{{t}_{N}}={{t}_{f}}$ where, for simplicity, we take
${{t}_{j}}-{{t}_{j-1}}=h$, for $j=1,...,N$.  A "hat" notation is
also used to distinguish the numerically estimated state and control
values from their actual counterparts that are "hatless", e.g.,
${{\hat{x}}_{j}}=\hat{x}({{t}_{j}})$,
$\,{{\hat{u}}_{0,j}}={{\hat{u}}_{0}}({{t}_{j}})$,
${{\hat{u}}_{1,j}}={{\hat{u}}_{1}}({{t}_{j}})$and
${{\hat{v}}_{j}}=\hat{v}({{t}_{j}})$.  The collocation method used
here assumes triangular basis functions for the state and piecewise
constant basis functions (derivatives of triangular functions) for
the controls. Specifically, the estimated state is given by
\begin{equation}
\hat{x}(t)=\sum\limits_{j=0}^{N}{{{{\hat{x}}}_{j}}{{\varphi
}_{j}}(t)}
\end{equation}
where the ${{\hat{x}}_{j}}$'s are to be determined and the
triangular basis functions are given by
\begin{equation}
{{\varphi }_{j}}(t)=\left\{ \begin{array}{*{35}{l}}
\frac{t-{{t}_{j-1}}}{h}, & {{t}_{j-1}}<t\le
{{t}_{j}}  \\
\frac{{{t}_{j+1}}-t}{h}, & {{t}_{j}}<t\le
{{t}_{j+1}}  \\
0, & \text{elsewhere}  \\
\end{array} \right.
\end{equation}
We note two points:  the method is not restricted to using triangular
basis functions and each of the ${{\varphi }_{j}}(t)$'s is a time shift
of the previous one.

As summarized in \cite{Zefran1996}, the theoretical approach for
computing the controls is to extend the state space with new state
variables, ${{x}_{ext}}\in {{R}^{m+1}}$, whose derivative are the
desired controls, $u(t)\in {{R}^{m}}$ and $v(t)\in [0,1]\subset R$,
to be computed.  However, our choice of triangular basis functions
for the states renders the control inputs piecewise constant and we
simply solve directly for these (constant) control values.
Specifically, the estimates of the control inputs are given by
\begin{equation}
\left[ \begin{matrix}
\hat{u}(t)  \\
\hat{v}(t)  \\
\end{matrix} \right]=\sum\limits_{j=1}^{N}{\left[ \begin{matrix}
{{{\hat{u}}}_{j}}  \\
{{{\hat{v}}}_{j}}  \\
\end{matrix} \right]{{\psi }_{j}}}(t)
\end{equation}
where the piecewise constant basis functions are given by
\begin{equation}
{{\psi }_{j}}(t)=\left\{ \begin{matrix}
1 & {{t}_{j-1}}<t\le {{t}_{j}}  \\
0 & \text{elsewhere}  \\
\end{matrix} \right.
\end{equation}
The essence of the midpoint rule in the collocation method is to
enforce the constraints at the midpoints of each interval
$[{{t}_{j-1}},{{t}_{j}}]$ for $j=1,...,N$.  There results the
discretized embedded state dynamics
\begin{equation}
\begin{array}{c}
{{{\hat{x}}}_{j}} = {{{\hat{x}}}_{j-1}}+h\cdot
(1-{{{\hat{v}}}_{j}})\cdot {{f}_{0}}\left(
\frac{{{{\hat{x}}}_{j-1}}+{{{\hat{x}}}_{j}}}{2},{{{\hat{u}}}_{0j}}
\right) \\
\ \ \ +h\cdot {{{\hat{v}}}_{j}}\cdot {{f}_{1}}\left(
\frac{{{{\hat{x}}}_{j-1}}+{{{\hat{x}}}_{j}}}{2},{{{\hat{u}}}_{1j}}
\right)
\end{array}
\end{equation}
for $j=1,...,N$, with ${{f}_{0}}(\cdot )$ and ${{f}_{1}}(\cdot )$
the discretized state dynamics in modes -0 and -1, respectively.
Thus, the solution to the EOCP is given by the following NLP:
Minimize
\begin{equation}
\begin{array}{c}
{{{\hat{J}}}_{E}} = {{C}_{bat}}(\centerdot ){{\left(
{{\overset{\hat{\ }}{\mathop{SOC}}\,}_{N}}-SO{{C}^{NOM}}
\right)}^{2}} \\
\ \ \ +\sum\limits_{j=1}^{N}{\frac{1}{2}h\left\{ {{L}_{E}}\left(
{{t}_{j}},{{{\hat{x}}}_{j}},{{{\hat{u}}}_{0j}},{{{\hat{u}}}_{1j}},{{{\hat{v}}}_{j}},{{{\hat{p}}}_{j}}
\right) \right.} \\
\ \ \ \left. +{{L}_{E}}\left(
{{t}_{j-1}},{{{\hat{x}}}_{j-1}},{{{\hat{u}}}_{0j}},{{{\hat{u}}}_{1j}},{{{\hat{v}}}_{j}},{{{\hat{p}}}_{j}}
\right) \right\}
\end{array}
\label{ZEqnNum822958}
\end{equation}
over the controls $({{\hat{u}}_{j}},{{\hat{v}}_{j}})\in \Omega
\times [0,1]$, subject to equation (33) and all other equality/power
flow constraints represented as
$g({{\hat{x}}_{j-1}},{{\hat{x}}_{j}},{{\hat{u}}_{j}},{{\hat{v}}_{j}},{{\hat{p}}_{j}})=0$.
Here ${{L}_{E}}(\cdot )$ is the integrand of equation properly
discretized and
${{\overset{\lower0.5em\hbox{$\smash{\scriptscriptstyle\frown}$}}{p}}_{j}}$
represents the various power flows in the model.

\subsection{Nonlinear Model Predictive Control}

The NMPC solution strategy in this study uses a moving four-second
predictive-window with the control applied over one-second
sub-interval.  In particular the NLP for $t={{t}_{j}}$ is solved over
$[{{t}_{j}},{{t}_{j+1}},{{t}_{j+2,}}{{t}_{j+3}},{{t}_{j+4}}]$ instead
of the entire driving cycle.  The resulting control at each iteration
is applied only over $[{{t}_{j}},{{t}_{j+1}}]$ to the system model.
Using the constant controls computed by NMPC algorithm, the system
model is then simulated over $[{{t}_{j}},{{t}_{j+1}}]$ to obtain an
updated state at ${{t}_{j+1}}$.  This updated state represents what
would be measured in a real-world implementation of the NMPC control
\cite{Martinsen2004, Nagy2004}.

We let ${{t}_{f,j}}$ denote the final time of each NMPC iteration,
and hence is the final time in the PI (\ref{ZEqnNum948398}).
Further for the NMPC strategy the coefficient penalizing the
deviation from nominal SOC, ${{C}_{bat}}({{t}_{f,j}})$ is linearly
interpolated according to the equation
\begin{equation}
{{C}_{bat}}({{t}_{f,j}})=({{{t}_{f,j}}}/{{{t}_{f}}}\;)C_{bat}^{nom}
\label{ZEqnNum911666}
\end{equation}
Otherwise, the NMPC control will try to maintain the SOC at
$SO{{C}^{NOM}}$ over each iteration unduly restricting the use of
battery power.

Difficulties arise for real world implementation since this strategy
presumes knowledge of ${{t}_{f}}$.  Nevertheless, it can be entered by
drivers, or becomes an adaptive function of the recent history of the
vehicle's power consumption.  The estimation problem of ${{t}_{f}}$ is
beyond the scope of this paper.

In the beginning of each NMPC window, we assume the knowledge of the
current road grade (e.g., through an accelerometer or future GPS).  The
control algorithms assume that over each partition of the NMPC window,
the road grade is constant at the value at the beginning of the NMPC
window.  Although the road grade may change over the NMPC window, since
the control is only applied over the first partition, after which a new
measurement is taken, potential error is believed negligible.

\section{SIMULATION RESULTS}

The NLP of the PMCP is applicable to various numerical solvers such as
AIMMS, TOMLAB, etc.  In this study, a sequential quadratic programming
(SQP) based NLP solver, \textit{fmincon}, is adopted.  This section
details simulation results for PHEV tracking different driving
profiles.  The first set of simulations compare the overall hybrid
optimal control and the NMPC tracking the sawtooth driving profile with
road grades.  The next simulation details the PHEV tracking the
standard 765 s EPA Highway driving profile to which is added a
sinusoidal road grade to better exercise the controller performance.
The second simulation looks at the 600 s US06 FTP supplemental driving
schedule.  The numerical solutions to various driving profiles are
shown below.

\subsection{Optimal and NMPC Tracking of Sawtooth Velocity
Profile with Road Grades: Cases 1 and 2}

The performances of the overall hybrid optimal control and the NMPC
strategies are compared using the sawtooth driving profile suggested
in \cite{Uthaichana2008}.  The sawtooth profile demands higher rates
of acceleration/deceleration than typical driving cycles. To further
test the limits of performance of the vehicle powertrain, positive
sinusoidal road grades are superimposed in this paper.

The coefficients of the PI (\ref{ZEqnNum948398}) are ${{C}_{V}}=10$,
${{C}_{ICE}}={{10}^{-3}}$, ${{C}_{FR}}={{10}^{-4}}$, and
${{C}_{bat}}(\centerdot )=C_{bat}^{nom}={{10}^{5}}$.  This study
compares the performances of two hybrid optimal control strategies.
Case 1 uses a control that minimizes the PI of (\ref{ZEqnNum948398})
with the above coefficients over the entire driving cycle and is
thus optimal over the driving cycle.  Case 2 constructs an NMPC
version that can be compared to the optimal solution.
\begin{figure}[h]
\centerline{\includegraphics[width=226pt]{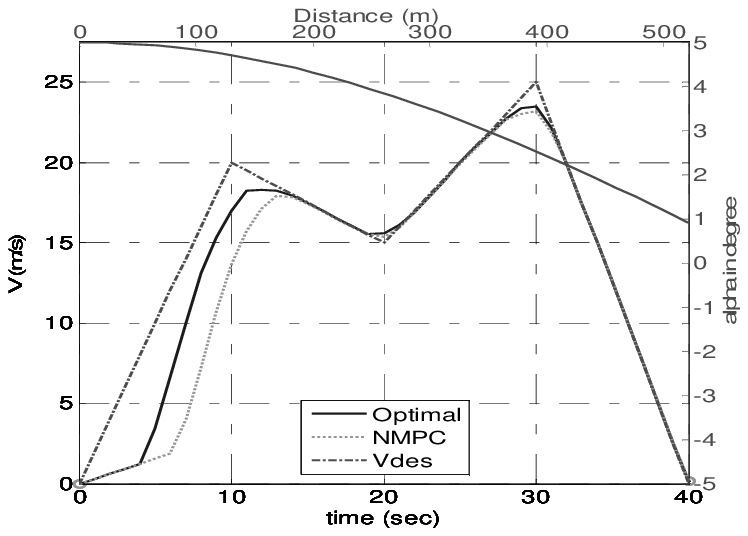}}
\caption{Velocity tracking performance for overall optimal hybrid
control and NMPC tracking sawtooth profile with road grades}
\label{Fig 5}
\end{figure}
Figure \ref{Fig 5} shows that initially the vehicle in both cases
fail to provide perfect tracking.  Initial tracking error is due to
the insufficient available propelling power from the ICE and ED both
of which already operate at their maximum levels, as shown in Fig.
\ref{Fig 6} for case 1.  As per Fig. \ref{Fig 6a}, the ICE is off at
startup (propelling from ED alone) due to the closed-loop local
control constraint, until a minimum operating speed of 800 RPM is
reached. Further, the tracking error during the first 14 s for the
NMPC version is larger due to the fact that the NMPC decides to turn
the ICE on slightly later.  The rest of the driving profiles except
the peak at 30 s can be tracked relatively well in both cases.

\begin{figure}[h!]
\begin{center}
\subfigure[ICE power]{
\includegraphics[width=226pt]{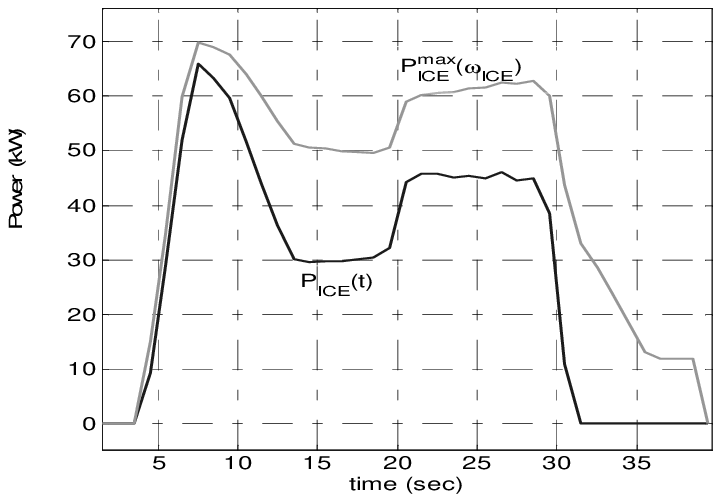}\label{Fig 6a}} \\
\subfigure[ED
power]{\includegraphics[width=226pt]{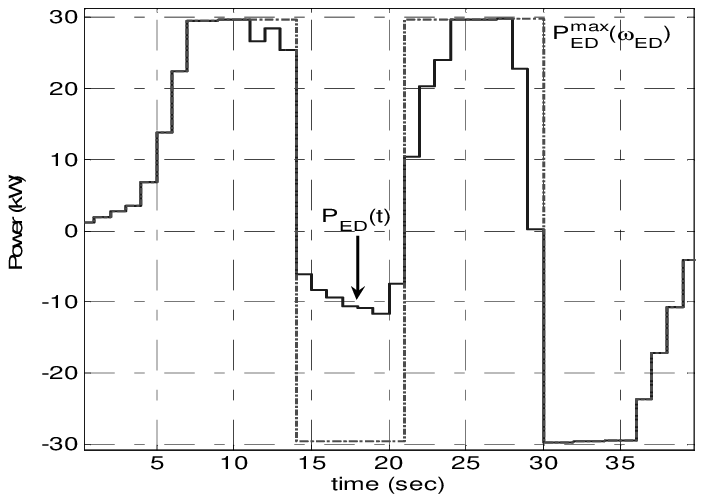}\label{Fig 6b}}
\end{center}
\caption{Output power for tracking sawtooth profile with road grades
in case 1} \label{Fig 6}
\end{figure}

\begin{figure}[ht!]
\begin{center}
\subfigure[ICE
power]{\includegraphics[width=226pt]{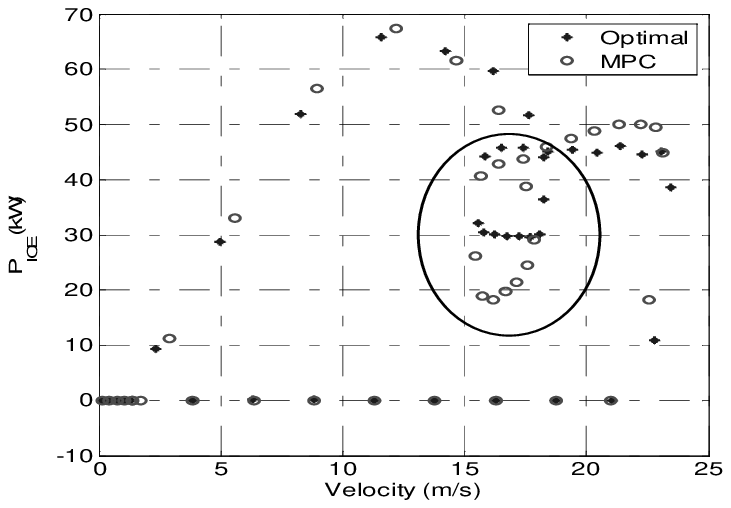}\label{Fig 7a}}
\subfigure[ED
power]{\includegraphics[width=226pt]{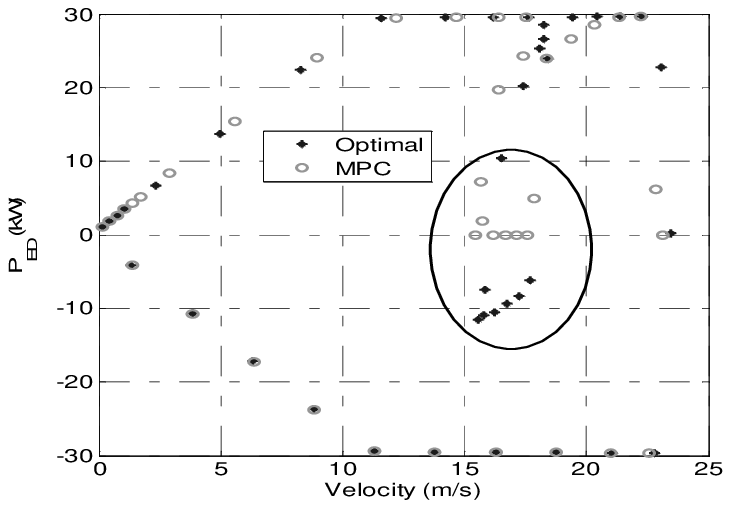}\label{Fig 7b}}
\end{center}
\label{Fig 7} \caption{Output profiles vs vehicle's velocity for
overall optimal hybrid control and NMPC tracking sawtooth with road
grades showing ICE power and ED power}
\end{figure}

After 30 s, when the desired velocity is decreasing and the road
grades are positive but small; the ICE is turned off (Fig. \ref{Fig
6a}) and the ED provides negative power to charge the battery in
both cases as shown in Fig. \ref{Fig 6b}. Between 14 and 21 s the
ICE in case 1 provides power both to charge the battery and to the
wheels, while in case 2, the ICE only provides propelling power. The
ICE profiles in this specific region for the NMPC are lower and
circled in Fig. \ref{Fig 7a}.  This contributes to slightly better
fuel economy for the NMPC version of 14.7 MPG vs. 14.2 MPG for case
1.

Simultaneously, the ED in case 1 operates as a generator (mode-1)
more often, specifically between 14 and 21 s, as depicted in Fig.
\ref{Fig 6b} and Fig. \ref{Fig 7b}.  The NMPC hands down a
non-generating decision in this time interval because of a
relatively lower penalty (coefficient) on the deviation from nominal
SOC according to (\ref{ZEqnNum911666}).

\begin{figure}[h!]
\centerline{\includegraphics[width=226pt]{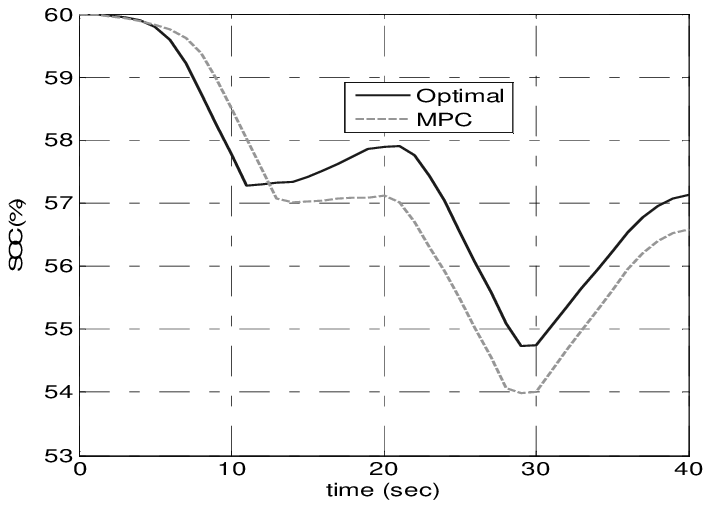}}
\caption{Battery SOC for overall optimal hybrid control and NMPC
tracking sawtooth with road grades} \label{Fig 8}
\end{figure}

Figure \ref{Fig 8} depicts the battery SOC profiles.  Case 1 shows
consistency with the mode of operation and returns to the vicinity
of the 60\% nominal level at the end of the cycle.

\begin{figure}[h!]
\centerline{\includegraphics[width=216pt]{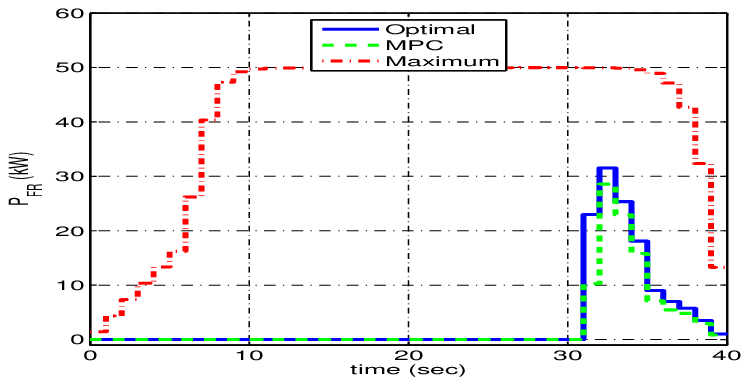}}
\caption{Frictional-braking power for overall optimal hybrid control
and NMPC tracking sawtooth with road grades} \label{Fig 9}
\end{figure}

In case 2, SOC deviation is only mildly penalized allowing lower
levels than in case 1 during the first half cycle.  In the second
half of the driving cycle, with an increasing penalty on the
deviation from nominal SOC, the NMPC tries to recharge the battery
and returns SOC to the level slightly lower than in case 1.

During the last 10 s, the ED operates as a generator in both cases
as expected to provide regenerative braking power.  Further, the
demanded power is more (negative) than the maximum level that the ED
can deliver in both cases as shown in Fig. \ref{Fig 6b}.  To achieve
the desired velocity tracking, the extra kinetic energy is expended
in frictional braking as shown in Fig. \ref{Fig 9}.  Indeed,
frictional braking is only used during the last 10 s of the driving
profile as one would expect.

\subsection{NMPC Tracking of EPA Highway Driving Cycle with Road
Grade: Case 3}

\begin{figure}[h!]
\centerline{\includegraphics[width=226pt]{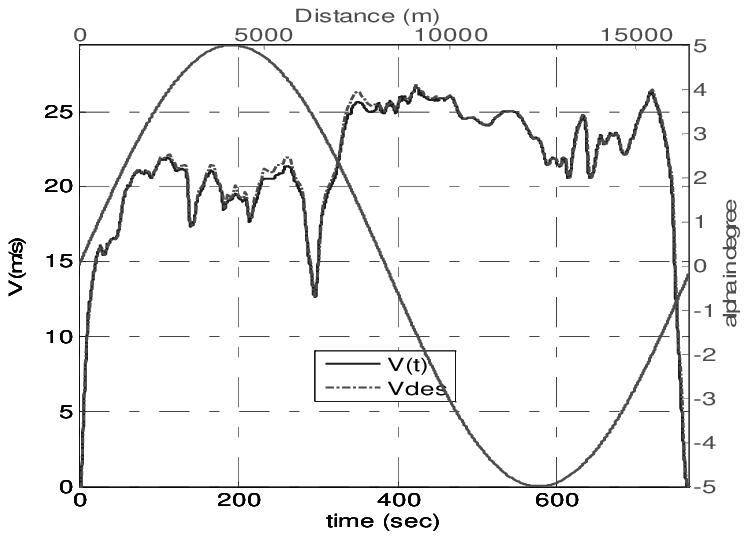}}
\caption{Velocity tracking performances for case 3 on a road grades
whose angle in degrees is indicated by the sinusoidal curve with
values in degrees on the right vertical axis} \label{Fig 10}
\end{figure}

In this case study, the vehicle is to follow the EPA highway driving
profile for 765 seconds on a nonzero road grades as shown in Fig.
\ref{Fig 10}. Also shown (solid sinusoidal) in Fig. \ref{Fig 10}, is
the road grades, which has a positive angle (uphill) over the first
382.5 s, and then has a negative angle (downhill) over the final
382.5 s

The coefficients of the PI (\ref{ZEqnNum948398}) are the same as
those used in the sawtooth driving profiles with a sliding penalty
on the deviation of the SOC from nominal as $t \to {{t}_{f}}=765$ s
Fig. 10 shows that the NMPC strategy provides nearly perfect
tracking for this case.

\begin{figure}[h!]
\begin{center}
\subfigure[Mode of
operation]{\includegraphics[width=226pt]{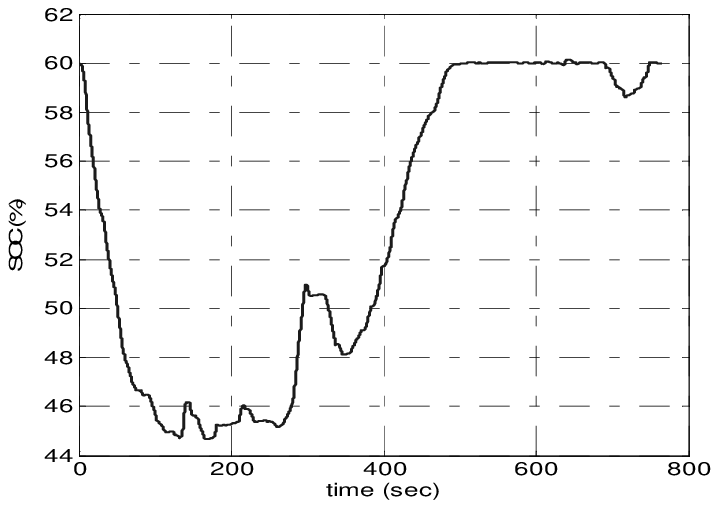}\label{Fig 11a}}
\subfigure[Battery
SOC]{\includegraphics[width=229pt]{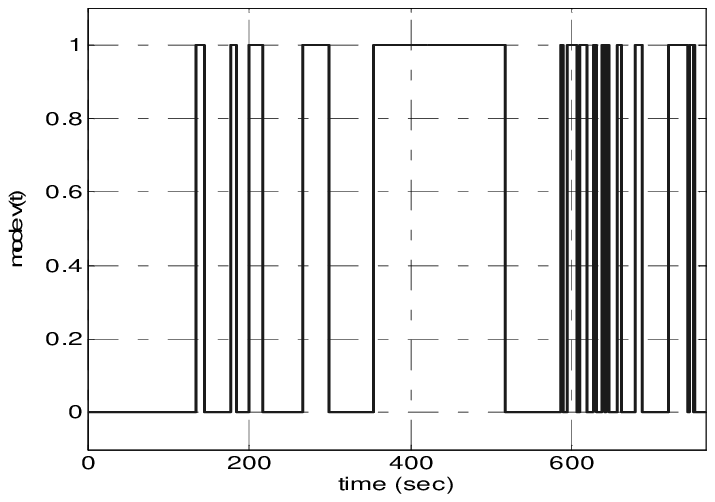}\label{Fig 11b}}
\end{center}
\label{Fig 11} \caption{NMPC strategy for tracking EPA highway
profile showing Mode of operation and Battery SOC Profiles}
\end{figure}

Since there is very little penalty on battery usage initially (see
(\ref{ZEqnNum911666})) and a relatively significant penalty on fuel
consumption, mode 0 is active for 135 s (Fig. \ref{Fig 11a}) at
which the desired velocity starts to decrease.

Initially, the ED supplies roughly 20 kW of propelling power, and
eventually decays to near zero at 135 s.  In contrast, the ICE is
initially off, and gradually ramps to supply propulsion demand reach
40 kW at roughly 135 s (Fig. \ref{Fig 13}).  In other words, the
NMPC strategy decides to drain the power from the ED-battery pack to
as low as 45\% during this time period as shown in Fig. \ref{Fig
11b}.
\begin{figure}[h!]
\centerline{\includegraphics[width=226pt]{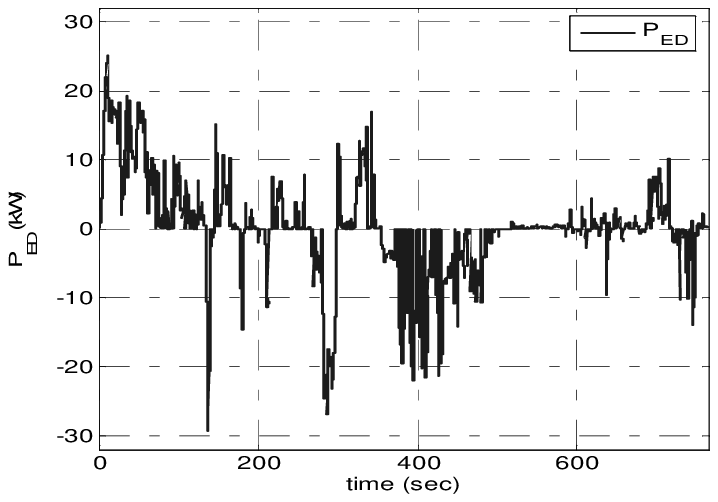}} \caption{ED
output power profiles for PHEV tracking EPA highway using NMPC
strategy} \label{Fig 12}
\end{figure}

After 60 seconds, the ICE power contribution continues to ramp and
then remains relatively constant at about 45 kW until 400 seconds as
per Fig. \ref{Fig 13}.  Over this period, fuel efficiency is
relatively high according to the ICE efficiency map in Fig. \ref{Fig
4}.
\begin{figure}[h!]
\centerline{\includegraphics[width=226pt]{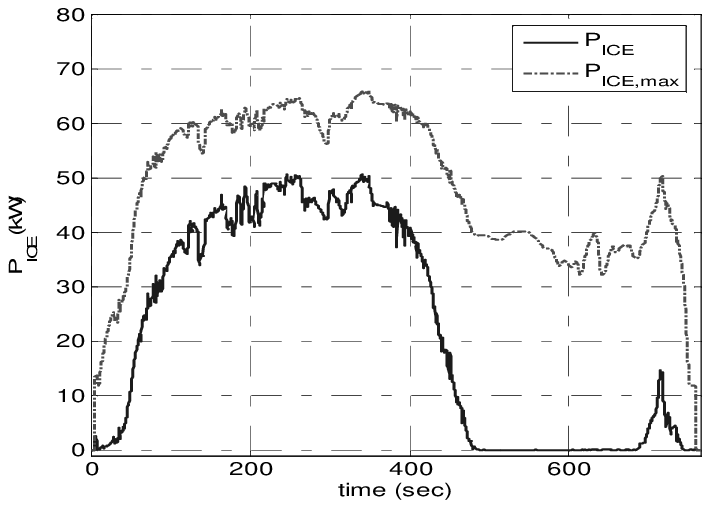}} \caption{ICE
power usages for PHEV tracking EPA highway using NMPC strategy}
\label{Fig 13}
\end{figure}

Around 260 sec, the road profile becomes less demanding and the
vehicle operates in the generating mode more often.  Further, after
382.5 sec, the road grades are negative, which leads to regenerative
braking.  Figure \ref{Fig 12} shows that the ED operates more
frequently in the generating mode for recharging the battery back to
its nominal value of 0.6.  Further, there is little need for
propelling power to maintain perfect tracking, hence the ICE is off
most of the time in the second half of the driving cycle.
\begin{figure}[h!]
\centerline{\includegraphics[width=226pt]{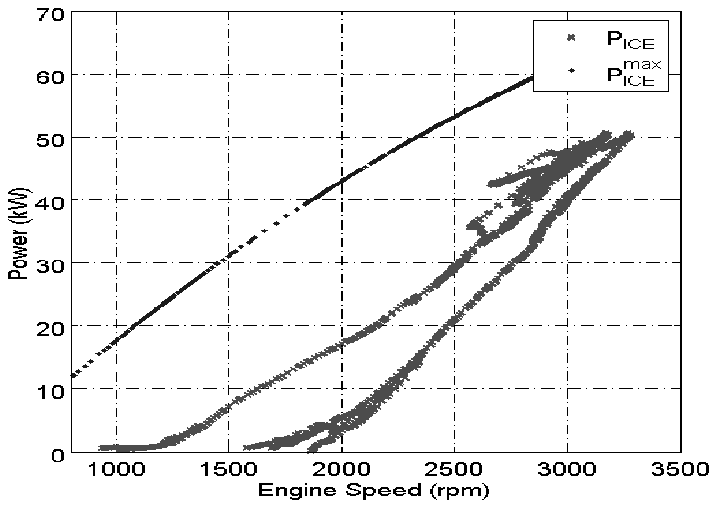}}
\caption{Trajectories of ICE power vs. engine speed for PHEV
tracking EPA highway using NMPC strategy showing
acceleration/deceleration hysteresis} \label{Fig 14}
\end{figure}

\begin{figure}[h!]
\centerline{\includegraphics[width=226pt]{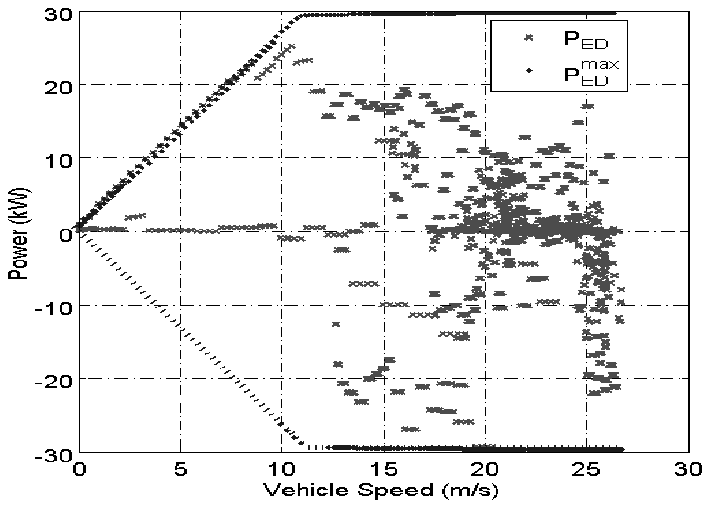}}
\caption{Trajectories of ED power vs. vehicle for PHEV tracking EPA
highway profile using NMPC} \label{Fig 15}
\end{figure}

The map of the ICE power over the engine speeds' range shows denser
data in the fuel efficient region.  Hysteresis appears as the result
of the dynamics (lag) in the engine power during the acceleration
and deceleration as shown in Fig. \ref{Fig 14}.  Figure \ref{Fig 15}
shows the operation of the ED over the range of vehicle speeds.

In the presence of the sinusoidal road grade over the cycle, the
fuel economy is 27.5 mpg, lower than the 32 mpg (similar to current
hybrid SUVs on the market) when the road is flat due to (i)
frictional braking losses during the negative road grades to
maintain velocity tracking and (ii) inefficiencies in energy
recovery from regenerative braking. The simulation result in the
case of vehicle tracking the EPA highway velocity profile on a flat
road is not included in this paper due to space limitations, but it
can be found in \cite{Uthaichana2006}.

One concludes that the NMPC strategy performs very well while
sustaining the constraint on the final SOC.  Further, the resulting
power distributions in this case study is used to justify the
concept of five-to-two mode reduction as shown in Fig. \ref{fig2}.
The 40 kW line is drawn in the figure to indicate a rough ICE power
level that is fuel efficient for a medium engine-speed range.

\subsection{MPC Tracking of the US06 supplemental FTP Driving
Profile: Case 4}

\begin{figure}[h!]
\centerline{\includegraphics[width=226pt]{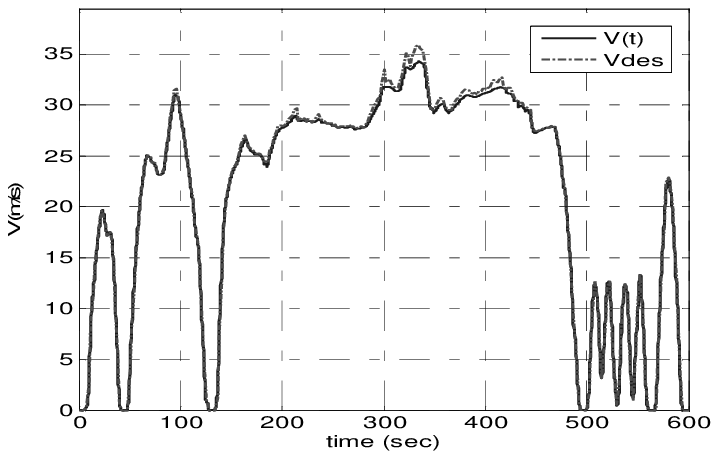}} \caption{US06
supplemental FTP driving profile and vehicle tracking performance}
\label{Fig 16}
\end{figure}

In this case study, we again use NMPC to track the 600 second US06
supplemental FTP driving schedule, which demands higher
accelerations and more aggressive velocity variation/limits than the
standard EPA city and highway schedules as shown in Fig. \ref{Fig
16}. The coefficients of the PI are the same as those in the
previous cases. Figure \ref{Fig 16} again shows that the NMPC almost
perfectly tracks the desired velocity profile.

\begin{figure}[h!]
\begin{center}
\subfigure[ED output
power]{\includegraphics[width=226pt]{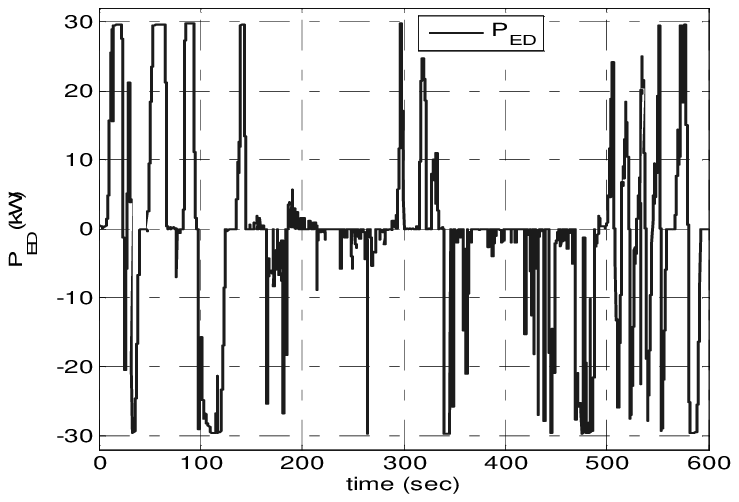}\label{Fig 17a}}
\subfigure[ICE output
power]{\includegraphics[width=226pt]{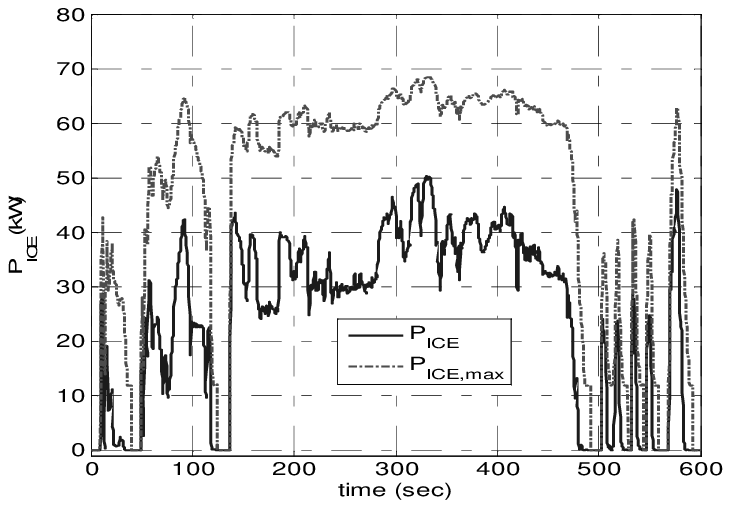}\label{Fig 17b}}
\end{center}
\caption{Plots of ED output power and ICE output power profiles in
case 4} \label{Fig 17}
\end{figure}

The aggressive nature of the velocity profile forces the ED to
operate close to maximum power levels in both modes (Fig. \ref{Fig
17a}) while the power usage of the ICE mimics (Fig. \ref{Fig 17b})
the shape of the velocity profile.

To meet the transient acceleration demands, NMPC puts the ED in the
motoring mode during acceleration and in the generating mode during
deceleration.  During the non-transient power demand, the middle of
the driving cycle, the ED is often off or provides relatively low
power and the ICE provides most of the power as it can operate
around its more fuel efficient level.  Overall, estimated fuel
efficiency averages at 23 mpg, a lower value than the EPA highway
profile of 32 due to the more aggressive power demands of this case
study.

\begin{figure}[h!]
\begin{center}
\subfigure[Battery
SOC]{\includegraphics[width=226pt]{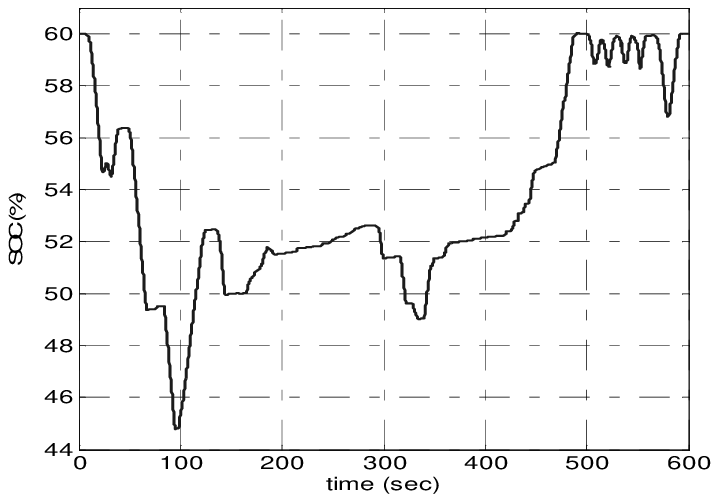}\label{Fig 18a}}
\subfigure[Mode of
operation]{\includegraphics[width=226pt]{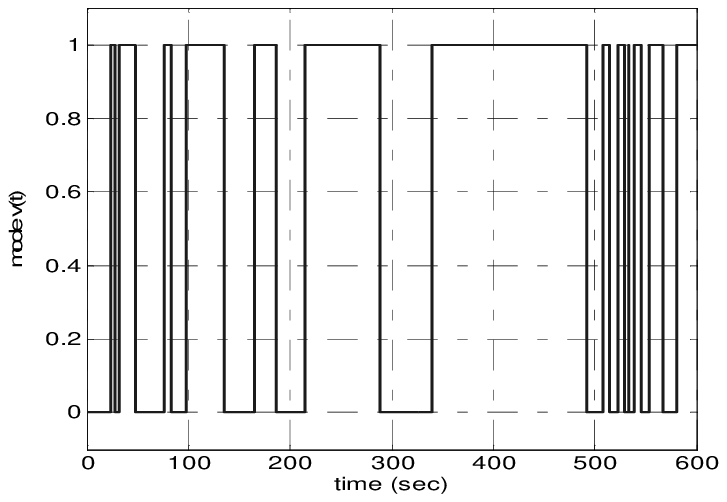}\label{Fig 18b}}
\end{center}
\caption{Case 4: Battery SOC and Mode of operation profiles}
\label{Fig 18}
\end{figure}

During the first 100 sec of the driving schedule, the HEV operates
primarily in the motoring mode draining the battery to an SOC of
about 45 percent as shown in Fig. \ref{Fig 18a}.  The aggressive
acceleration driving profile in concert with the relatively low
penalty on the SOC deviation from the nominal level contributes to
the NMPC selection of this strategy.  Afterward, as the penalty on
the SOC deviation increases linearly toward the final time of 600
sec, the SOC rises toward the nominal level of 0.6 with more
frequent time in mode 1.

Mode switching in both cases, using the strategy outlined in section
4.3, is reasonable, and consistent with the velocity variations in
the driving profile as shown in Fig. \ref{Fig 18b}.

\begin{figure}[h!]
\centerline{\includegraphics[width=226pt]{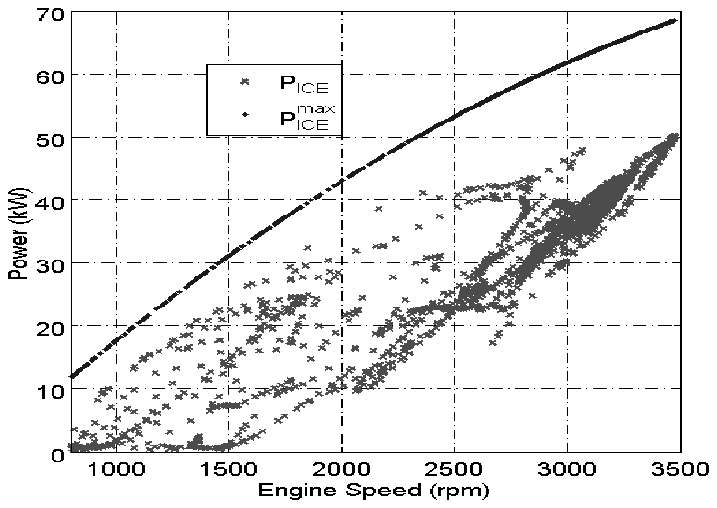}}
\caption{Trajectories of ICE power vs. engine speed for PHEV
tracking US06 FTP profile using NMPC strategy} \label{Fig 19}
\end{figure}

\begin{figure}[h!]
\centerline{\includegraphics[width=226pt]{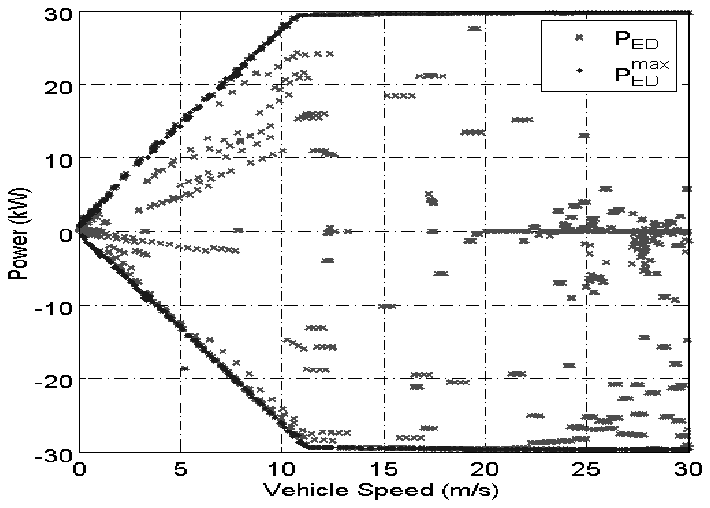}}
\caption{Trajectories of ED power vs. vehicle for PHEV tracking US06
FTP profile using NMPC strategy} \label{Fig 20}
\end{figure}

Similar to case 3, the nature of the modified speed envelope
strategy for the CVT dictates higher engine speed for higher ICE
power as depicted in Fig. \ref{Fig 19}.  Dense engine data is
located around 40 kW.  In contrast to case 3, wider spread in the
trajectories of the ICE on the power-vs-speed map can be observed in
this case. Frequent operations on the ED maximum power-vs-speed
envelope (in both motoring and generating) are also observed in this
case, as depicted in Fig. \ref{Fig 20}. These characteristics in the
ICE and the ED power profiles reflect greater variety of demands in
the US06 FTP profile.

\subsection{CONCLUSIONS}

An application of the hybrid optimal control for solving the power
management control problem (PMCP) in a parallel electric hybrid
vehicle (PHEV) has been illustrated in this paper.  The advantages
outlined in \cite{Bengea2005} motivate solving for the
optimal/suboptimal power flow as the solution to the embedded
version of the problem, i.e., the EOCP.  The solution of the EOCP
can be obtained via a number of numerical techniques.

In this study, the numerical solution is obtained by converting the
original infinite dimensional problem into a finite dimensional
nonlinear programming problem (NLP) using the direct collocation
technique.  Then, the resulting NLP is solved via a sequential
quadratic programming algorithm.  Requiring a short predictive window,
the NMPC strategy is applied to solve the PMCP for the sawtooth, EPA
highway, and US06 supplemental FTP driving profiles.  The simulations
show that the NMPC can track the driving profiles quite well unless
there is insufficient available power to achieve the tracking as
illustrated in case of the sawtooth profile with positive road grades.
In practice, vehicle control strategies are often PID and maps/look-up
tables based.  The resulting NMPC profiles can help providing
additional information on how to improve the existing look-up tables or
tuning the gain-scheduling maps in the PID based controllers.

\bibliographystyle{ieeetr}
\bibliography{hevbib_exported}

\begin{thebibliography}{10}

\bibitem{Erdinc2009}
O.~Erdinc, B.~Vural, and M.~Uzunoglu, ``A wavelet-fuzzy logic based energy
  management strategy for a fuel cell/battery/ultra-capacitor hybrid vehicular
  power system,'' {\em Journal of Power Sources}, vol.~In Press, Corrected
  Proof, 2009.

\bibitem{Pesaran2009}
A.~Pesaran, J.~Gonder, and M.~Keyser, ``Ultracapacitor applications and
  evaluation for hybrid electric vehicles,'' in {\em Advanced Capacitors World
  Summit 2009}, (La Jolla, California), 2009.

\bibitem{Barbir2005}
F.~Barbir, H.~Gorgun, and X.~Wang, ``Relationship between pressure drop and
  cell resistance as a diagnostic tool for pem fuel cells,'' {\em Journal of
  Power Sources}, vol.~141, no.~1, pp.~96--101, 2005.

\bibitem{Cao2009}
J.~B. Cao and B.~G. Cao, ``Neural network sliding mode control based on on-line
  identification for electric vehicle with ultracapacitor-battery hybrid
  power,'' {\em International Journal of Control, Automation and Systems},
  vol.~7, no.~3, pp.~409--418, 2009.

\bibitem{Feroldi2009}
D.~Feroldi, M.~Serra, and J.~Riera, ``Energy management strategies based on
  efficiency map for fuel cell hybrid vehicles,'' {\em Journal of Power
  Sources}, vol.~190, no.~2, pp.~387--401, 2009.

\bibitem{Shen2008}
Q.~Shen, M.~Hou, X.~Yan, D.~Liang, Z.~Zang, L.~Hao, Z.~Shao, Z.~Hou, P.~Ming,
  and B.~Yi, ``The voltage characteristics of proton exchange membrane fuel
  cell (pemfc) under steady and transient states,'' {\em Journal of Power
  Sources}, vol.~179, no.~1, pp.~292--296, 2008.

\bibitem{Springer1991}
T.~E. Springer, T.~A. Zawodzinski, and S.~Gottesfeld, ``Polymer electrolyte
  fuel cell model,'' {\em Journal of the Electrochemical Society}, vol.~138,
  no.~8, pp.~2334--2342, 1991.

\bibitem{Ao2007}
G.~Q. Ao, J.~A. Qiang, H.~Zhong, L.~Yang, and B.~Zhuo, ``Exploring the fuel
  economy potential of isg hybrid electric vehicles through dynamic
  programming,'' {\em International Journal of Automotive Technology}, vol.~8,
  pp.~781--790, Dec 2007.

\bibitem{Lin2003}
C.-C. Lin, P.~Huei, J.~W. Grizzle, and K.~Jun-Mo, ``Power management strategy
  for a parallel hybrid electric truck,'' {\em IEEE Trans. Control Syst.
  Technol.}, vol.~11, no.~6, pp.~839--49, 2003.

\bibitem{Pisu2007}
P.~Pisu and G.~Rizzoni, ``A comparative study of supervisory control strategies
  for hybrid electric vehicles,'' {\em IEEE Transactions on Control Systems
  Technology}, vol.~15, no.~3, pp.~506--518, 2007.

\bibitem{Rizzoni1999}
G.~Rizzoni, L.~Guzzella, and B.~M. Baumann, ``Unified modeling of hybrid
  electric vehicle drivetrains,'' {\em IEEE/ASME Transactions on Mechatronics},
  vol.~4, no.~3, pp.~246--257, 1999.

\bibitem{Scordia2005}
J.~Scordia, M.~Desbois-Renaudin, R.~Trigui, B.~Jeanneret, F.~Badin, and
  C.~Plasse, ``Global optimisation of energy management laws in hybrid vehicles
  using dynamic programming,'' {\em International Journal of Vehicle Design},
  vol.~39, no.~4, pp.~349--367, 2005.

\bibitem{Lee2008}
J.~M. Lee and J.~H. Lee, ``Value function-based approach to the scheduling of
  multiple controllers,'' {\em Journal of Process Control}, vol.~18, no.~6,
  pp.~533--542, 2008.

\bibitem{Powell2009}
W.~B. Powell, ``What you should know about approximate dynamic programming,''
  {\em Naval Research Logistics}, vol.~56, no.~3, pp.~239--249, 2009.

\bibitem{Kim2010}
N.~Kim, S.~Cha, and H.~Peng, ``Optimal control of hybrid electric vehicles
  based on pontryagin's minimum principle,'' {\em IEEE Transactions on Control
  Systems Technology}, 2010.
\newblock cited By (since 1996) 0; Article in Press.

\bibitem{Schouten2003}
N.~J. Schouten, M.~A. Salman, and N.~A. Kheir, ``Energy management strategies
  for parallel hybrid vehicles using fuzzy logic,'' {\em Control Engineering
  Practice}, vol.~11, no.~2, pp.~171--177, 2003.

\bibitem{Xiong2009}
W.~Xiong, Y.~Zhang, and C.~Yin, ``Optimal energy management for a
  series-parallel hybrid electric bus,'' {\em Energy Conversion and
  Management}, vol.~50, no.~7, pp.~1730--1738, 2009.

\bibitem{Kermani2008}
S.~Kermani, S.~Delprat, R.~Trigui, and T.~M. Guerra, ``Predictive energy
  management of hybrid vehicle,'' in {\em IEEE Vehicle Power and Propulsion
  Conference}, 2008 IEEE Vehicle Power and Propulsion Conference, VPPC 2008,
  (Harbin), 2008.

\bibitem{Bemporad1999}
A.~Bemporad and M.~Morari, ``Control of systems integrating logic, dynamics,
  and constraints,'' {\em Automatica}, vol.~35, no.~3, pp.~407--27, 1999.

\bibitem{Bengea2005}
S.~C. Bengea and R.~A. DeCarlo, ``Optimal control of switching systems,'' {\em
  Automatica}, vol.~41, no.~1, pp.~11--27, 2005.

\bibitem{Neuman1973}
C.~P. Neuman and A.~Sen, ``A suboptimal control algorithm for constrained
  problems using cublic splines,'' {\em Automatica}, vol.~9, pp.~601--13, 1973.

\bibitem{Schaefer2007}
A.~Schäfer, P.~Kühl, M.~Diehl, J.~Schlöder, and H.~G. Bock, ``Fast reduced
  multiple shooting methods for nonlinear model predictive control,'' {\em
  Chemical Engineering and Processing: Process Intensification}, vol.~46,
  no.~11, pp.~1200--1214, 2007.

\bibitem{VonStryk1993}
O.~Von~Stryk, ``Numerical solution of optimal control problems by direct
  collocation. in optimal control—calculus of variations, optimal control
  theory, and numerical methods,'' {\em International Series of Numerical
  Mathematics}, vol.~111, pp.~129--143, 1993.

\bibitem{Zefran1996}
M.~{\v Z}efran, {\em Continuous Methods for Motion Planning}.
\newblock PhD thesis, University of Pennsylvania, Philadelphia, 1996.

\bibitem{Uthaichana2005}
K.~Uthaichana, S.~Bengea, and R.~DeCarlo, ``Suboptimal supervisory level power
  flow control of a hybrid electric vehicle,'' in {\em Proc. of IFAC World
  Congress}, (Prague), July 2005.

\bibitem{Diehl2006}
M.~Diehl, H.~G. Bock, H.~Diedam, and P.~B. Wieber, ``Fast direct multiple
  shooting algorithms for optimal robot control,'' in {\em Fast Motions in
  Biomechanics and Robotics} (M.~Diehl and K.~Mombaur, eds.), vol.~340 of {\em
  Lecture Notes in Control and Information Sciences}, pp.~65--93, Springer,
  2006.

\bibitem{DeHaan2006}
D.~DeHaan and M.~Guay, ``A new real-time perspective on non-linear model
  predictive control,'' {\em Journal of Process Control}, vol.~16, no.~6,
  pp.~615--624, 2006.

\bibitem{Martinsen2004}
F.~Martinsen, L.~T. Biegler, and B.~A. Foss, ``A new optimization algorithm
  with application to nonlinear mpc,'' {\em Journal of Process Control},
  vol.~14, no.~8, pp.~853--865, 2004.

\bibitem{Floudas1995}
C.~A. Floudas, {\em Nonlinear and Mixed-Integer Optimization: Fundamentals and
  Applications}.
\newblock New York: Oxford University Press, 1995.

\bibitem{Nemhauser1988}
G.~L. Nemhauser and L.~A. Wolsey, {\em Integer and Combinatorial Optimization}.
\newblock Discrete Mathematics and Optimiazation, New York: Wiley-Interscience,
  1988.

\bibitem{Uthaichana2006}
K.~Uthaichana, {\em Modeling and control of a parallel hybrid electric
  vehicle}.
\newblock PhD thesis, Purdue, West Lafayette, December 2006 2006.

\bibitem{Uthaichana2008}
K.~Uthaichana, S.~Bengea, R.~DeCarlo, S.~Pekarek, and M.~{\v Z}efran, ``Hybrid
  model predictive control tracking of a sawtooth driving profile for an hev,''
  in {\em Proceedings of the American Control Conference}, 2008 American
  Control Conference, ACC, (Seattle, WA), pp.~967--974, 2008.

\bibitem{Gill1981}
P.~Gill, W.~Murray, and M.~Wright, {\em Practical Optimization}.
\newblock Academic Press, 1981.

\bibitem{Pfiffner2001}
R.~Pfiffner and L.~Guzzella, ``Optimal operation of cvt-based powertrains,''
  {\em Int. J. Robust Nonlinear Control}, vol.~11, no.~11, pp.~1003--1021,
  2001.

\bibitem{Coleman2008}
M.~Coleman, W.~G. Hurley, and C.~K. Lee, ``An improved battery characterization
  method using a two-pulse load test,'' {\em IEEE Transactions on Energy
  Conversion}, vol.~23, no.~2, pp.~708--713, 2008.

\bibitem{Rao2003}
R.~Rao, S.~Virudhula, and D.~Rakhmatov, ``Battery models for energy aware
  system design,'' {\em IEEE Computer}, vol.~36, pp.~1019--1030, 2003.

\bibitem{Agarwal2009}
V.~Agarwal, K.~Uthaichana, R.~DeCarlo, and L.~H. Tsoukalas, ``Development and
  validation of a battery model useful for discharging and charging power
  control and lifetime estimation,'' {\em IEEE Transactions on Energy
  Conversion}, vol.~25, pp.~821--835, August 2010.

\bibitem{Pekarek2005}
S.~Pekarek, K.~Uthaichana, S.~Bengea, R.~DeCarlo, and M.~{\v Z}efran,
  ``Modeling of an electric drive for a hev supervisory level power flow
  control problem,'' in {\em 2005 IEEE Vehicle Power and Propulsion
  Conference}, (IEEE Cat. No.05EX1117C), (Chicago, IL, USA), IEEE, 2005.

\bibitem{Anatone2005}
M.~Anatone, R.~Cipollone, and A.~Sciarretta, ``Control-oriented modeling and
  fuel optimal control of a series hybrid bus,'' in {\em SAE Technical Papers
  2005-01-1163}, 2005.

\bibitem{Jeon2002}
S.-I. Jeon, S.-T. Jo, Y.-I. Park, and J.-M. Lee, ``Multi-mode driving control
  of a parallel hybrid electric vehicle using driving pattern recognition,''
  {\em Journal of Dynamic Systems, Measurement and Control, Transactions of the
  ASME}, vol.~124, no.~1, pp.~141--149, 2002.

\bibitem{Koot2005}
M.~Koot, J.~T. B.~A. Kessels, B.~de~Jager, W.~P. M.~H. Heemels, P.~P.~J.
  van~den Bosch, and M.~Steinbuch, ``Energy management strategies for vehicular
  electric power systems,'' {\em IEEE Transactions on Vehicular Technology},
  vol.~54, no.~3, pp.~771--82, 2005.

\bibitem{Heywood1988}
J.~Heywood, {\em Internal Combustion Engine Fundamentals}.
\newblock McGraw-Hill, 1988.

\bibitem{Giua2001}
A.~Giua, C.~Seatzu, and C.~Van Der~Mee, ``Optimal control of switched
  autonomous linear systems,'' in {\em Proceedings of the IEEE Conference on
  Decision and Control}, vol.~3, (Orlando, FL), pp.~2472--77, Institute of
  Electrical and Electronics Engineers Inc., 2001.

\bibitem{Hedlund1999}
S.~Hedlund and A.~Rantzer, ``Optimal control of hybrid systems,'' in {\em
  Proceedings of the IEEE Conference on Decision and Control}, vol.~4,
  (Phoenix, AZ, USA), pp.~3972--77, IEEE, Piscataway, NJ, USA, 1999.

\bibitem{Wei2007}
S.~Wei, K.~Uthaichana, M.~{\v Z}efran, R.~A. DeCarlo, and S.~Bengea,
  ``Applications of numerical optimal control to nonlinear hybrid systems,''
  {\em Nonlinear Analysis: Hybrid Systems}, vol.~1, no.~2, pp.~264--279, 2007.

\bibitem{Nagy2004}
Z.~K. Nagy and F.~Allg\"ower, ``Nonlinear model predictive control: From
  chemical industry to microelectronics,'' in {\em Proceedings of the IEEE
  Conference on Decision and Control}, vol.~4 of {\em 2004 43rd IEEE Conference
  on Decision and Control (CDC)}, (Nassau), pp.~4249--4254, 2004.

\end{thebibliography}

\end{document}